\documentclass[12pt]{article}
\usepackage{latexsym, amssymb}

\usepackage[pdftex]{graphicx}

\DeclareMathAlphabet\gothic{U}{euf}{m}{n}

\makeatletter
\def\eqnarray{\stepcounter{equation}\let\@currentlabel=\theequation
\global\@eqnswtrue
\tabskip\@centering\let\\=\@eqncr
$$\halign to \displaywidth\bgroup\hfil\global\@eqcnt\z@
  $\displaystyle\tabskip\z@{##}$&\global\@eqcnt\@ne
  \hfil$\displaystyle{{}##{}}$\hfil
  &\global\@eqcnt\tw@ $\displaystyle{##}$\hfil
  \tabskip\@centering&\llap{##}\tabskip\z@\cr}

\def\endeqnarray{\@@eqncr\egroup
      \global\advance\c@equation\m@ne$$\global\@ignoretrue}

\def\@yeqncr{\@ifnextchar [{\@xeqncr}{\@xeqncr[5pt]}}
\makeatother

\begin{document}
\bibliographystyle{tom}

\newtheorem{lemma}{Lemma}[section]
\newtheorem{thm}[lemma]{Theorem}
\newtheorem{cor}[lemma]{Corollary}
\newtheorem{voorb}[lemma]{Example}
\newtheorem{rem}[lemma]{Remark}
\newtheorem{prop}[lemma]{Proposition}
\newtheorem{stat}[lemma]{{\hspace{-5pt}}}
\newtheorem{obs}[lemma]{Observation}
\newtheorem{defin}[lemma]{Definition}

\newenvironment{remarkn}{\begin{rem} \rm}{\end{rem}}
\newenvironment{exam}{\begin{voorb} \rm}{\end{voorb}}
\newenvironment{defn}{\begin{defin} \rm}{\end{defin}}
\newenvironment{obsn}{\begin{obs} \rm}{\end{obs}}

\newenvironment{emphit}{\begin{itemize} }{\end{itemize}}

\newcommand{\gota}{\gothic{a}}
\newcommand{\gotb}{\gothic{b}}
\newcommand{\gotc}{\gothic{c}}
\newcommand{\gote}{\gothic{e}}
\newcommand{\gotf}{\gothic{f}}
\newcommand{\gotg}{\gothic{g}}
\newcommand{\gothh}{\gothic{h}}
\newcommand{\gotk}{\gothic{k}}
\newcommand{\gotm}{\gothic{m}}
\newcommand{\gotn}{\gothic{n}}
\newcommand{\gotp}{\gothic{p}}
\newcommand{\gotq}{\gothic{q}}
\newcommand{\gotr}{\gothic{r}}
\newcommand{\gots}{\gothic{s}}
\newcommand{\gotu}{\gothic{u}}
\newcommand{\gotv}{\gothic{v}}
\newcommand{\gotw}{\gothic{w}}
\newcommand{\gotz}{\gothic{z}}
\newcommand{\gotA}{\gothic{A}}
\newcommand{\gotB}{\gothic{B}}
\newcommand{\gotG}{\gothic{G}}
\newcommand{\gotL}{\gothic{L}}
\newcommand{\gotS}{\gothic{S}}
\newcommand{\gotT}{\gothic{T}}

\newcommand{\mn}{\marginpar{\hspace{1cm}*} }
\newcommand{\mnn}{\marginpar{\hspace{1cm}**} }

\newcommand{\mnq}{\marginpar{\hspace{1cm}*???} }
\newcommand{\mnnq}{\marginpar{\hspace{1cm}**???} }

\newcounter{teller}
\renewcommand{\theteller}{\Roman{teller}}
\newenvironment{tabel}{\begin{list}%
{\rm \bf \Roman{teller}.\hfill}{\usecounter{teller} \leftmargin=1.1cm
\labelwidth=1.1cm \labelsep=0cm \parsep=0cm}
                      }{\end{list}}

\newcounter{tellerr}
\renewcommand{\thetellerr}{(\roman{tellerr})}
\newenvironment{subtabel}{\begin{list}%
{\rm  (\roman{tellerr})\hfill}{\usecounter{tellerr} \leftmargin=1.1cm
\labelwidth=1.1cm \labelsep=0cm \parsep=0cm}
                         }{\end{list}}
\newenvironment{ssubtabel}{\begin{list}%
{\rm  (\roman{tellerr})\hfill}{\usecounter{tellerr} \leftmargin=1.1cm
\labelwidth=1.1cm \labelsep=0cm \parsep=0cm \topsep=1.5mm}
                         }{\end{list}}

\newcommand{\Ni}{{\bf N}}
\newcommand{\Ri}{{\bf R}}
\newcommand{\Ci}{{\bf C}}
\newcommand{\Ti}{{\bf T}}
\newcommand{\Zi}{{\bf Z}}
\newcommand{\Fi}{{\bf F}}

\newcommand{\proof}{\mbox{\bf Proof} \hspace{5pt}} 
\newcommand{\remark}{\mbox{\bf Remark} \hspace{5pt}}
\newcommand{\ruimte}{\vskip10.0pt plus 4.0pt minus 6.0pt}

\newcommand{\simh}{{\stackrel{{\rm cap}}{\sim}}}
\newcommand{\ad}{{\mathop{\rm ad}}}
\newcommand{\Ad}{{\mathop{\rm Ad}}}
\newcommand{\Aut}{\mathop{\rm Aut}}
\newcommand{\arccot}{\mathop{\rm arccot}}
\newcommand{\capp}{{\mathop{\rm cap}}}
\newcommand{\rcapp}{{\mathop{\rm rcap}}}
\newcommand{\diam}{\mathop{\rm diam}}
\newcommand{\divv}{\mathop{\rm div}}
\newcommand{\dist}{\mathop{\rm dist}}
\newcommand{\codim}{\mathop{\rm codim}}
\newcommand{\RRe}{\mathop{\rm Re}}
\newcommand{\IIm}{\mathop{\rm Im}}
\newcommand{\Tr}{{\mathop{\rm Tr}}}
\newcommand{\Vol}{{\mathop{\rm Vol}}}
\newcommand{\card}{{\mathop{\rm card}}}
\newcommand{\supp}{\mathop{\rm supp}}
\newcommand{\sgn}{\mathop{\rm sgn}}
\newcommand{\essinf}{\mathop{\rm ess\,inf}}
\newcommand{\esssup}{\mathop{\rm ess\,sup}}
\newcommand{\Int}{\mathop{\rm Int}}
\newcommand{\Leibniz}{\mathop{\rm Leibniz}}
\newcommand{\lcm}{\mathop{\rm lcm}}
\newcommand{\loc}{{\rm loc}}

\newcommand{\mod}{\mathop{\rm mod}}
\newcommand{\spann}{\mathop{\rm span}}
\newcommand{\one}{1\hspace{-4.5pt}1}

\newcommand{\DWR}{}

\hyphenation{groups}
\hyphenation{unitary}

\newcommand{\tfrac}[2]{{\textstyle \frac{#1}{#2}}}

\newcommand{\ca}{{\cal A}}
\newcommand{\cb}{{\cal B}}
\newcommand{\cc}{{\cal C}}
\newcommand{\cd}{{\cal D}}
\newcommand{\ce}{{\cal E}}
\newcommand{\cf}{{\cal F}}
\newcommand{\ch}{{\cal H}}
\newcommand{\ci}{{\cal I}}
\newcommand{\ck}{{\cal K}}
\newcommand{\cl}{{\cal L}}
\newcommand{\cm}{{\cal M}}
\newcommand{\cn}{{\cal N}}
\newcommand{\co}{{\cal O}}
\newcommand{\cs}{{\cal S}}
\newcommand{\ct}{{\cal T}}
\newcommand{\cx}{{\cal X}}
\newcommand{\cy}{{\cal Y}}
\newcommand{\cz}{{\cal Z}}

\newcommand{\wtozp}{W^{1,2}\raisebox{10pt}[0pt][0pt]{\makebox[0pt]{\hspace{-34pt}$\scriptstyle\circ$}}}
\newlength{\hightcharacter}
\newlength{\widthcharacter}
\newcommand{\covsup}[1]{\settowidth{\widthcharacter}{$#1$}\addtolength{\widthcharacter}{-0.15em}\settoheight{\hightcharacter}{$#1$}\addtolength{\hightcharacter}{0.1ex}#1\raisebox{\hightcharacter}[0pt][0pt]{\makebox[0pt]{\hspace{-\widthcharacter}$\scriptstyle\circ$}}}
\newcommand{\cov}[1]{\settowidth{\widthcharacter}{$#1$}\addtolength{\widthcharacter}{-0.15em}\settoheight{\hightcharacter}{$#1$}\addtolength{\hightcharacter}{0.1ex}#1\raisebox{\hightcharacter}{\makebox[0pt]{\hspace{-\widthcharacter}$\scriptstyle\circ$}}}
\newcommand{\scov}[1]{\settowidth{\widthcharacter}{$#1$}\addtolength{\widthcharacter}{-0.15em}\settoheight{\hightcharacter}{$#1$}\addtolength{\hightcharacter}{0.1ex}#1\raisebox{0.7\hightcharacter}{\makebox[0pt]{\hspace{-\widthcharacter}$\scriptstyle\circ$}}}

 \thispagestyle{empty}

\vspace*{1.5cm}
\begin{center}
{\Large{\bf Uniqueness of diffusion on domains}}\\[2mm] 
{\Large{\bf  with rough boundaries }}  \\[2mm]
\large  Juha Lehrb{\"a}ck$^1$  and Derek W. Robinson$^2$ \\[2mm]

\normalsize{31st March 2015}
\end{center}

\vspace{5mm}

\begin{center}
{\bf Abstract}
\end{center}

\begin{list}{}{\leftmargin=1.5cm \rightmargin=1.5cm \listparindent=10mm 
   \parsep=0pt}
   \item
Let $\Omega$ be a  domain in~$\Ri^d$
and $h(\varphi)=\sum^d_{k,l=1}(\partial_k\varphi, c_{kl}\partial_l\varphi)$ a quadratic form on $L_2(\Omega)$ with
domain $C_c^\infty(\Omega)$ where
 the  $c_{kl}$ are real symmetric $L_\infty(\Omega)$-functions with $C(x)=(c_{kl}(x))>0$ for almost all  $x\in \Omega$.
Further assume there are $a, \delta>0$ such that  $a^{-1}d_\Gamma^{\;\delta}\,I\leq C\leq a\,d_\Gamma^{\;\delta}\,I$ for $d_\Gamma\leq 1$ 
where $d_\Gamma$ is the Euclidean distance to the boundary $\Gamma$ of $\Omega$.

We assume that   $\Gamma$ is Ahlfors $s$-regular  and if $s$, the Hausdorff dimension of $\Gamma$,  is larger or equal to $d-1$ 
we also assume  a mild uniformity property for $\Omega$ in the neighbourhood of one $z\in\Gamma$.
Then we establish that
$h$ is Markov unique, i.e.\ it has a unique Dirichlet form extension, if and only if $\delta\geq 1+(s-(d-1))$. 
The result applies to forms on   Lipschitz domains or on  a wide class of domains with $\Gamma$ a self-similar fractal.
In particular it applies to 
the interior or exterior of the von Koch snowflake curve in $\Ri^2$
or the complement of a uniformly disconnected set in $\Ri^d$.

\vfill
\end{list}

\vspace{1cm}


\vspace{5mm}

\noindent
\begin{tabular}{@{}cl@{\hspace{10mm}}cl}
1. &  Department of Mathematics 
 & 
 2. & Centre for Mathematics  \\
&\hspace{10mm} and Statistics  & 
  & \hspace{5mm} and its Applications 
   \\
& University  of Jyvaskyla
& 
  &Mathematical Sciences Institute 
  \\
&  PO Box 35 (MaD) & 
 & Australian National University  \\
& FI-400014 University of Jyvaskyla & {}
 & Canberra, ACT 0200 
 \\
& Finland & {}
  & Australia
  \\[5pt]
& juha.lehrback@jyu.fi & {}
 &derek.robinson@anu.edu.au \\
\end{tabular}

\newpage
\setcounter{page}{1}

\section{Introduction}\label{S1}

The theory of diffusion processes has a distinct probabilistic character and is most naturally
studied on $L_1$-spaces.
Consequently much of the analysis of such processes has relied on methods of stochastic differential equations
or stochastic integration.
Our aim, however, is to examine symmetric diffusion problems on domains of Euclidean space with the
techniques of functional analysis and semigroup theory.
In particular we focus on the  characterization of uniqueness of the $L_1$-theory
on domains with rough or fragmented boundaries.
First we formulate the problem of diffusion as a problem of finding extensions of a given 
elliptic operator which generate  semigroups with the general characteristics suited to the description of diffusion.

Let $\Omega$ be a domain in $ \Ri^d $, i.e.\ a non-empty open connected subset,  with boundary $\partial\Omega$
and $S=\{S_t\}_{t\geq0}$ a strongly continuous, positive, contraction semigroup on $L_1(\Omega)$.
If the positive normalized functions in $L_1(\Omega)$ are viewed as probability distributions then $S$ has the basic properties
required for description of their evolution with time.
For brevity we refer to $S$ as a diffusion semigroup.
We define $S$ to be symmetric if 
\begin{equation}
(S_t\varphi, \psi)=(\varphi, S_t\psi)
\label{euni1.0}
\end{equation}
for all $\varphi\in L_1(\Omega)$, all $\psi\in L_1(\Omega)\cap L_\infty(\Omega)$ and all $t\geq0$.
It follows that $S$ extends by continuity from $ L_1(\Omega)\cap L_\infty(\Omega)$ to a weakly$^*$ continuous semigroup on $L_\infty(\Omega)$
which we also denote by $S$.
The extended semigroup is automatically equal to the adjoint semigroup $S^*=\{S^*_t\}_{t\geq0}$.
Then $S$ can be defined  on $L_p(\Omega)$ for each $p\in\langle1,\infty\rangle$ by interpolation.
In particular $S$ is a self-adjoint, positive, contraction semigroup on $L_2(\Omega)$.
If $H$ is the positive, self-adjoint generator of $S$ it  then follows from the Beurling--Deny criteria (see, for example, \cite{RS4}) that the corresponding
quadratic form $h(\varphi)=\|H\varphi\|_2^2$ with $\varphi\in D(h)=D(H^{1/2})$ is a Dirichlet form.
Therefore the semigroup $S$ is submarkovian, i.e.\ if $0\leq \varphi\leq \one$ then  $0\leq S_t\varphi\leq \one$ for all $t>0$, by the   theory of Dirichlet forms 
\cite{BH} \cite{FOT}.

Next define the operator  $H_0$ on the domain $D(H_0)=C_c^\infty(\Omega)$ by 
\begin{equation}
H_0\varphi
= -\sum^d_{k,l=1}\partial_k\, c_{kl} \, \partial_l\varphi
\label{euni1.10}
\end{equation}
where $c_{kl}=c_{lk}\in W^{1,\infty}(\Omega)$ are 
real and  the matrix of coefficients $C(x)=(c_{kl}(x))>0$ for all $x\in\Omega$ in the sense of matrix order.
The corresponding diffusion problem consists of classifying all extensions of $H_0$ to $L_1(\Omega)$  which generate symmetric diffusion semigroups.
One can establish the existence of at least one such extension by quadratic form techniques.
Let $h_0$ be the positive, quadratic, form associated with $H_0$ on $L_2(\Omega)$, i.e.\ 
\begin{equation}
h_0(\varphi)
=(\varphi, H_0\varphi)= \sum^d_{k,l=1}(\partial_k\varphi, c_{kl} \, \partial_l\varphi)
\label{euni1.1}
\end{equation}
for all $\varphi\in D(h_0)= C_c^\infty(\Omega)$.
Since $H_0$ is a symmetric operator on $L_2(\Omega)$ the form $h_0$ is  closable  and the closure, 
which we denote by $h_D$, is automatically a Dirichlet form \cite{BH} \cite{FOT}.
The corresponding positive, self-adjoint operator $H_D$, the Friedrichs' extension of $H_0$,  generates a   positive, contraction semigroup $S^D$
on $L_2(\Omega)$ which extends to a similar semigroup on each of the $L_p$-spaces.
The extension to $L_1(\Omega)$ automatically satisfies the symmetry relation (\ref{euni1.0}).
Therefore  $H_D$  generates a symmetric diffusion  semigroup on $L_1(\Omega)$.
The extension $H_D$  corresponds to Dirichlet boundary conditions on $\partial\Omega$.
But the same argument establishes that  each Dirichlet form extension of $h_0$ determines the generator of a symmetric diffusion semigroup on $L_1(\Omega)$.
Therefore there is a one-to-one correspondence between extensions of $H_0$ on $L_1(\Omega)$ which generate symmetric diffusion semigroups 
and Dirichlet form extensions of $h_0$ on $L_2(\Omega)$.
The classification of extensions of $H_0$ which generate symmetric diffusion semigroups on $L_1(\Omega)$ is now reduced 
to the more amenable and transparent  problem of  classifying the Dirichlet form extensions of $h_0$ on~$L_2(\Omega)$.

The Dirichlet form extensions of $h_0$ have a fundamental ordering property.
The closure $h_D$ is the smallest Dirichlet form extension of $h_0$ but there is also a largest such extension $h_N$.
The maximal extension $h_N$ is defined on the domain 
\[
D(h_N)=\{\varphi\in W^{1,2}_{\rm loc}(\Omega):\,\Gamma(\varphi)+\varphi^2\in L_1(\Omega)\}\;,
\]
where $\Gamma(\varphi)=\sum^d_{k,l=1}c_{kl}(\partial_k\varphi)(\partial_l\varphi)$ is the {\it carr\'e du champ}, by  setting 
\[
h_N(\varphi)=\int_\Omega dx\,\Gamma(\varphi)(x)
\]
for  $\varphi\in D(h_N)$.
Then $h_N$ is a Dirichlet form and the associated operator $H_N$ is the extension of $H_0$ corresponding to generalized Neumann  boundary conditions.
But  if $k$ is a general Dirichlet form extension of $h_0$ then $D(h_D)\subseteq D(k)\subseteq D(h_N)$ (see \cite{FOT}, Section~3.3,
\cite{RSi4}, Theorem~1.1, or \cite{RSi5}, Theorem~2.1).
Thus $h_N\leq k\leq h_D$ in the sense of ordering of quadratic forms.
Clearly all the Dirichlet form extensions $k$ are equal in the interior of $\Omega$  and differ only by their behaviour at the boundary.
If $d=1$ then a classification of the possible extensions can be extracted from the general analyis of Feller \cite{Feller2}.
But a classification of the extensions in terms of boundary conditions seems well beyond reach if $d\geq 2$.
The multi-dimensional problem is  complicated by the wide range of geometric possibilities
for   $\Omega$ and the wide variety of possible boundary conditions.
Nevertheless these observations give a direct approach to the characterization  of 
uniqueness of a Dirichlet form extension  and consequently the uniqueness of the solution to the diffusion problem.

First define the form  $h_0$ to be Markov unique if the closure $h_D$ is the  unique Dirichlet form extension.
Thus $h_0$ is Markov unique if and only if $h_D=h_N$.
It was established in  \cite{RSi4} \cite{RSi5} (see also Section~\ref{S2}) that this latter condition is  equivalent to 
 the boundary $\partial\Omega$ having capacity zero  measured with respect to the form $h_N$.
This criterion is  a property  which depends
on the degeneracy of the coefficients $c_{kl}$  near the boundary
together with the regularity and uniformity properties of $\partial\Omega$. 
It does not, however, depend on any smoothness of the coefficients. 
Therefore in  the subsequent analysis of the uniqueness problem we replace the assumption $c_{kl}\in W^{1,\infty}(\Omega)$ by the weaker assumption $c_{kl}\in L_\infty(\Omega)$. 
Specifically we now assume that $h_0$ is defined by (\ref{euni1.1})  with real coefficients $c_{kl}=c_{lk}\in L_\infty(\Omega)$ and with $C(x)=(c_{kl}(x))>0$ for almost  all $x\in\Omega$.
Then for each compact subset $K\subset \Omega$ there is a $c_K>0$ such that $C(x)\geq c_K I$ for almost all $x\in K$.
Hence $h_0$ is closable (see \cite{MR}, Section~II.2b) and one can again define $h_D$ as the closure.
Moreover $h_N$, defined as above, is again a Dirichlet form (see \cite{OuR}, Proposition~2.1).
Therefore one can analyze the Markov uniqueness condition $h_D=h_N$ in this broader framework.
This equality depends critically on the behaviour of  the coefficients  on the boundary and we next 
formulate an appropriate degeneracy condition.

Let $d(x\,;y) $ denote the Euclidean distance from $x$ to $y$ 
and $B(x\,;r)$ the open  Euclidean ball with centre $x$ and radius $r$.
Further set $d_A(x)=\inf_{z\in A}d(x\,;z)$ for each non-empty subset $A\subset \Ri^d$.
Then, denoting the boundary of $\Omega$ by $\Gamma$, we
assume there is a  $\delta\geq 0$ and   for each bounded  non-empty subset $A\subset \Gamma$
there are $a,b,r>0$ such that 
\begin{equation}
a\,d_{\Gamma}(x)^{\delta}I\leq C(x)\leq b\,d_{\Gamma}(x)^{\delta}I
\label{euni1.30}
\end{equation}
for  almost all $x\in A_r$ where $A_r=\{x\in\overline\Omega: d_A(x)<r\}$ 
is  the (inner) $r$-neighbourhood of~$A$.

Next we place some mild geometric restraints on the domain $\Omega$.

First we suppose  that $\Gamma$ satisfies a property of  Ahlfors $s$-regularity.
Specifically we assume that there is a regular Borel measure $\mu$ on $\Gamma$ and an $ s>0$  such that 
for each subset  $A=\Gamma \cap B(x_0\,;R_0)$, with  $x_0\in\Gamma$ and $R_0>0$, one can choose
$c>0$ so that 
\begin{equation} 
c^{-1}\,r^s \leq \mu(A\cap B(x\,;r)) \leq c\, r^s
\label{euni1.3}
\end{equation}
for all $x\in A$ and $r\in \langle0, 2R_0\rangle$.
This is a locally uniform version of the Ahlfors  regularity  property used in the theory of metric spaces (see, for example, the monographs  \cite{DaS},  
\cite{Semm}, \cite{Hei} or \cite{MaT}). 
It implies that $\mu$ and the Hausdorff measure $\ch^s$ on $\Gamma$  are locally equivalent and $s=d_H(\Gamma)$, the Hausdorff dimension of~$\Gamma$.
The terminology regular  is somewhat misleading as an Ahlfors regular boundary can be quite `rough',
e.g.\ the  boundary of  the von Koch snowflake (see Section~\ref{S5}) is Ahlfors regular.
Condition (\ref{euni1.3}) does,  however,  imply that $\Gamma$ is regular in the sense that each of the subsets $\Gamma_{\!x,r}=\Gamma \cap B(x\,;r)$ with $x\in\Gamma$ has  Hausdorff dimension $s$.

Secondly, if $s\geq d-1$  we  assume a  local form  of the uniformity property  introduced by 
Martio and Sarvas \cite{MarS}.
If $z\in\Gamma$ and $R>0$ set $\Omega_{z,R}=\Omega\cap B(z\,;R)$.
Then $\Omega_{z,R}$ is defined to be $\Omega$-uniform 
 if there is a  $\sigma\geq1$ such that for all 
$x,y\in\Omega_{z,R}$ there is a rectifiable  curve
$\gamma\colon[0,1]\mapsto \Omega$ with $\gamma(0)=x$, $\gamma(1)=y$ and  length at most  $\sigma \,d(x\,;y)$
such that $d_\Gamma(w)\geq \sigma^{-1}(d(x\,;w)\wedge d(w\,;y))$ for all  $w\in\gamma([0,1])$.
Note that the curve $\gamma$ is in $\Omega$ but is not constrained to $\Omega_{z,R}$.

The local uniformity condition has two elements. 
First,  if an arbitrary pair of points $x,y\in\Omega_{z,R}$ can be joined by a rectifiable curve  in $\Omega$
then $\Omega_{z,R}$ must 
 belong to
a connected component of $\Omega\cap B(z,\sigma R)$
(but note that $\Omega_{z,R}$ need not be connected).
Secondly, it is necessary for  the detailed properties of the curves~$\gamma$ that the boundary subset $\Gamma_{\!z,R}=\Gamma\cap B(z\,;R)$
has the characteristics of the boundary of a uniform domain.
For example, if $d=2$ then   outward pointing parabolic cusps, inward pointing antennae or slits which separate locally are all forbidden.

The foregoing assumptions allow a rather simple characterization of Markov uniqueness in terms of the order of degeneracy $\delta$ of the coefficients of the form $h_0$
at the boundary $\Gamma$ and the Hausdorff dimension $s$.
 
 \begin{thm}\label{tuni1.1}
 Let $\Omega$ be a  domain in $ \Ri^d $ with boundary $\Gamma$.
Assume $\Gamma$  satisfies the Ahlfors $s$-regularity property $(\ref{euni1.3})$ with $s\in \langle0, d\rangle$.
Further,
 if $s\in[d-1,d\rangle$ assume  there is a $z\in\Gamma$  and an $R>0$ such that  $\Omega_{z,R}$ is  $\Omega$-uniform. 
Finally assume   the coefficients of the form $h_0$ satisfy the degeneracy condition $(\ref{euni1.30})$
 for $\delta\geq 0$.
 
\smallskip

 Then the  form $h_0$ is Markov unique, i.e.\ $h_D=h_N$, if and only if 
  $\delta\geq 1+(s-(d-1))$.
 \end{thm}

 Theorem~\ref{tuni1.1} is a straightforward illustration of our principal results.
 In the sequel (see Section~\ref{S4}) we describe situations with the index of regularity $s$ and the order of degeneracy $\delta$ 
 taking different values on distinct components and faces of the boundary.
 This introduces a number of extra complications but the basic elements of the proofs are already contained in the proof of the simpler theorem.

Despite the relative simplicity of Theorem~\ref{tuni1.1} it does cover a variety of interesting examples.
First if $\Omega$ is a Lipschitz domain then the regularity and uniformity assumptions of the theorem are valid with $s=d-1$ (see Section~\ref{S5}).
Therefore   $h_0$ is Markov unique if and only if $\delta\geq 1$.
Secondly, the theorem also applies to a  broad class of domains whose boundaries are self-similar fractals.
In particular it is applicable if $d=2$ and $\Omega$ is the interior, or exterior, of the von Koch snowflake.
Therefore $h_0$ is Markov unique if and only if $\delta\geq s$  with $s=\log4/\log3$, the Hausdorff dimension of the snowflake. 
 Thirdly, the conclusions of the theorem are stable under the subtraction of Ahlfors $s'$-regular subsets of the interior of $\Omega$ with $s'\leq s$
 (see Corollary~\ref{cuni4.1}).
Finally let   $\Gamma$ be  a uniformly disconnected subset of $\Ri^d$ (see, for example, \cite{Hei} Section~14.24).
MacManus, \cite{McM} page~275, observed that by the compactness argument   of V{\"a}is{\"a}l{\"a},  \cite{Vai} Theorem~3.6, the complement
 $\Omega=\Ri^d\backslash\Gamma$ is a uniform domain.
 (We give an explicit proof of this result in  Lemma~\ref{luni3.1}.)
Therefore  Theorem~\ref{tuni1.1} in combination with this observation   immediately gives the following corollary. 
  
  \begin{cor}\label{cuni1.1}
 Let  $\Omega=\Ri^d\backslash\Gamma$ where 
$\Gamma$  is a closed uniformly disconnected set satisfying the Ahlfors $s$-regularity property $(\ref{euni1.3})$  with $s\in \langle0, d\rangle$.
Further assume   the coefficients of the form $h_0$ satisfy the degeneracy condition $(\ref{euni1.30})$
 for $\delta\geq 0$.
 
\smallskip

 Then the  form $h_0$ is Markov unique, i.e.\ $h_D=h_N$, if and only if 
  $\delta\geq 1+(s-(d-1))$.
 \end{cor}

In particular if $\Gamma$ is the usual Cantor dust with granular ratio $\lambda\in\langle0,1/2\rangle$ then  $\Gamma$ is uniformly disconnected and $s=d_H(\Gamma)= d\log2/\log(1/\lambda)$. 
This is of interest as the Hausdorff dimension can take all values between $0$ and $d$ as $\lambda$ varies from $0$ to $1/2$.
Note that by setting $\delta=0$ and $C(x)=I$ one deduces that the  Laplacian  defined on $C_c^\infty(\Omega)$ is Markov unique if and only if $d_H(\Gamma)\leq d-2$.

\section{Preliminaries}
\label{S2}

In this section we gather some preliminary results which are needed in the proof Theorem~\ref{tuni1.1}.
First we recall some earlier results which characterize Markov uniqueness by  a zero capacity condition on the boundary $\Gamma$.
Secondly, we establish some implications of Ahlfors regularity and local uniformity of the  boundary.

The Markov uniqueness criterion $h_D=h_N$ is by definition equivalent to the density, with respect to the $D(h_N)$-graph norm 
$\psi\mapsto \|\psi\|_{D(h_N)}= (h_N(\psi)+\|\psi\|_2^2)^{1/2}$, of $C_c^\infty(\Omega)$ in $D(h_N)$.
But this criterion is a boundary condition and in earlier papers \cite{RSi4} \cite{RSi5} \cite{Rob8} it was established that it is equivalent to 
 $\Gamma$ having zero capacity relative to the form~$h_0$.
These earlier results were stated for forms with coefficients $c_{kl}\in W^{1, \infty}(\Omega)$ or $W^{1,\infty}_{\rm loc}(\Omega)$ but 
in fact no smoothness of the coefficients is necessary.
The property of importance is the density in $D(h_N)$, equipped with the graph norm,  of the subspace $(D(h_N)\cap L_\infty(\Omega))_c$ of  bounded functions in $D(h_N)$ with compact support in $\overline\Omega$.
This density property follows from the boundedness of the coefficients $c_{kl}$ of $h_0$.
In fact the density holds for large classes of coefficients which grow at infinity but fails in general (see \cite{Maz}, Section~2.7, or \cite{OuR}, Lemma~2.3, and \cite{Rob8}, Section~4).

The capacity of  a  subset $A$ of $\overline\Omega$ relative to the form $h_0$ is defined  by
\begin{eqnarray*}
\capp_{h_0}(A)=\inf\Big\{\;\|\psi\|_{D(h_N)}^2&&\;: \;\psi\in D(h_N), \;0\leq\psi\leq 1  \mbox{ and  there exists   an open set  }\nonumber\\[-5pt]
&& U\subset \Ri^d
\mbox{ such that } U\supseteq A
\mbox{ and } \psi=1 \mbox{ on } U\cap\Omega\;\Big\}
\;.
\end{eqnarray*}
The definition of $\capp_{h_0}(A)$ is analogous  to the  canonical definition  of the capacity associated with a Dirichlet form
\cite{BH}  \cite{FOT}
and  if $\Omega=\Ri^d$ the two definitions coincide.
The two capacities share many general characteristics.
If $A$ is measurable then $\capp_{h_0}(A)\geq |A|$  where $|A|$ denotes the Lebesgue measure of $A$. 
Moreover, the map $A\to\capp_{h_0}(A)$ is monotonic, if $A_n$ is an increasing family of measurable sets 
then $\capp_{h_0}(\cup_{n\geq 1}A_n)=\lim_{n\to\infty}\capp_{h_0}(A_n)$ and   if $A_n$ is a decreasing family of compact sets
then $\capp_{h_0}(\cap_{n\geq 1}A_n)=\lim_{n\to\infty}\capp_{h_0}(A_n)$. 

\smallskip

The following proposition is a slight extension of Theorem~1.2 of \cite{RSi4}.
\begin{prop}\label{puni2.1}
Let $\Omega$ be a  domain in $\Ri^d$ with boundary $\Gamma$ and $h_0$ the quadratic form with $L_\infty$-coefficients defined by $(\ref{euni1.1})$.

The following conditions are equivalent:
\begin{tabel}
\item\label{p1}
$h_D=h_N$,
\item\label{p2}
$\capp_{h_0}(\Gamma)=0$.
\end{tabel}
\end{prop}

Theorem~1.2 of \cite{RSi4} gives a similar statement for  $c_{kl}\in W^{1,\infty}(\Omega)$.
This smoothness property ensures that $h_0$ is the form of a symmetric operator.
But the argument used to establish the statement does not depend on smoothness. 
It is a quadratic form argument which applies equally well for $c_{kl}\in L_\infty(\Omega)$ 

The degeneracy conditions (\ref{euni1.30})  and regularity conditions (\ref{euni1.3})   have not been assumed in Proposition~\ref{puni2.1}.
The upper bounds of (\ref{euni1.30})  and the lower bounds of  (\ref{euni1.3})  are, however, critical for the subsequent verification of  the   criterion  $\capp_{h_0}(\Gamma)=0$ (see Section~\ref{S3}).
Another  crucial factor is the growth in volume of inner neighbourhoods $A_r=\{x\in\overline\Omega: d_A(x)<r\}$ of subsets $A$ of $\Gamma$.
The simplest estimates  on the growth are given by the following  lemma.

\begin{lemma}\label{luni2.2}
Let $A$ be a bounded non-empty subset of the boundary $\Gamma$ of the domain $\Omega$ and $\mu$ a regular  Borel measure on $A$.
\begin{tabel}
\item\label{luni2.2-1}
If there exist $a, R>0$ and $s\in\langle0, d\rangle$ such that 
$\mu(A\cap B(x\,;r))\geq a\,r^s$
for all $x\in A$ and all $r\in\langle0,R\rangle$
then there exists a $c>0$  such that 
\begin{equation}
| A_r|\leq c\, r^{d-s}
\label{euni2.2}
\end{equation}
for all $r\in\langle0, R\rangle$.
\item\label{luni2.2-2}
If  there exist 
$b, b', R'>0$ and $s\in\langle0,d\rangle$ such that 
$\mu(A\cap B(x\,;r))\leq b\,r^s$
and in addition $|\Omega\cap B(x\,; r)|\geq b'\,r^d$
for all $x\in A$ and all $r\in\langle0,R'\rangle$
then there exists a $c'>0$  such that 
\begin{equation}
|A_r|\geq  c'\,r^{d-s}
\label{euni2.3}
\end{equation}
for all $r\in\langle0,R'\rangle$.
\end{tabel}
\end{lemma}
\proof\  If $\Omega=\Ri^d$ and $A$  is a general bounded non-empty subset then the proposition follows from Lemma~2.1 in \cite{Sall}.
The proof of the latter lemma is based on a standard packing/covering argument.

If $\Omega\subset \Ri^d$  then Statement~\ref{luni2.2-1} is an immediate corollary of the $\Omega=\Ri^d$ case
because  $|A_r|\leq |\{x\in\Ri^d :d_A(x)\leq r\}|$.
The proof of Statement~\ref{luni2.2-2} is also a  corollary of Salli's argument but the  
lower bound on $|\Omega\cap B(x\,; r)|$  is essential.\hfill$\Box$

\bigskip

Since the assumptions of Lemma~\ref{luni2.2}.\ref{luni2.2-1}  imply that $|\{x\in\Ri^d:d_A(x)<r\}|\leq c\,r^{d-s}$ for all small $r>0$ it follows
in the limit $r\to0$ that $|A|=0$.
In particular if the boundary $\Gamma$ of the domain $\Omega$ is
 Ahlfors $s$-regular with $s\in\langle0,d\rangle$ then $|\Gamma|=0$.

The   lower bounds on $|A_r|$ given by the second statement  of  Lemma~\ref{luni2.2} depend on the  bounds $|\Omega\cap B(x\,; r)|\geq b'\,r^d$. 
These latter bounds are not generally valid but require some additional assumptions.
One general result in this direction is the following.

\begin{prop}\label{puni2.2}
Assume the boundary $\Gamma$ satisfies the Ahlfors $s$-regularity condition $(\ref{euni1.3})$  with  $s\in\langle0,d-1\rangle$.
If $A=\Gamma_{\!z,R}$ for some $z\in\Gamma$ and some $R>0$  then 
 there are  $c', R'>0$ such that 
\begin{equation}
|A_r|\geq  c'\,r^{d-s}
\label{euni2.40}
\end{equation}
 for all $r\in\langle0,R'\rangle$.
\end{prop}
\proof\
The proof follows by establishing that   $| \Omega\cap B(x\,;r)|=| \Ri^d\cap B(x\,;r)|$ for all $x\in A$ and $r>0$.

Set $B=B(x\,;r)$.
 If $\Omega'$ denotes the interior of the complement $\Omega^{\rm c}$ of $\Omega$ then $\Ri^d=\Omega\cup \Omega'\cup\Gamma$
and $|\Ri^d\cap B|=|\Omega\cap B|+|\Omega'\cap B|+|\Gamma\cap B|$.
But it follows from the Ahlfors regularity that $| \Gamma\cap B|=0$.
Now assume $\Omega'\cap B\neq \emptyset$.
Then $\Gamma\cap B$ separates $\Omega\cap B$ and $\Omega'\cap B$, i.e.\ each rectifiable curve in $B$ starting at $x_0\in \Omega\cap B$
and ending at $x_1\in  \Omega'\cap B$ has an intermediate point in $ \Gamma\cap B$.
Therefore $\Gamma\cap B$ has topological dimension $d-1$.
Consequently, the  Hausdorff dimension $s$ of $\Gamma\cap B$  is greater or equal to $d-1$ (see \cite{HuW}, Section~VII.4, or \cite{Hei}, Section~8.13).
But this is a contradiction so one must have $\Omega'\cap B=\emptyset$.
Hence $|\Omega\cap B|=|\Ri^d\cap B|= |B(0\,;1)|\,r^d$.
Now the lower bound (\ref{euni2.40})   follows from Statement~\ref{luni2.2-2} of Proposition~\ref{luni2.2}.
\hfill$\Box$

\bigskip

The lower bounds of Proposition~\ref{puni2.2} are independent of any uniformity property  of $\Omega$.
But if there is a $z\in\Gamma$ and $R>0$ such that $\Omega_{z,R}$ is $\Omega$-uniform then one has similar bounds for $s\in[d-1,d\rangle$  in the neighbourhood of the point of uniformity.

\begin{prop}\label{puni2.3}
Assume the boundary $\Gamma$ satisfies the Ahlfors $s$-regularity condition $(\ref{euni1.3})$ with  $s\in [d-1,d\rangle$.
Further  assume there is a $z\in\Gamma$  and an $R>0$ such  that $\Omega_{z,R}$  is $\Omega$-uniform.

It follows that  if  $A=\Gamma_{\!z,R/2}$ then there exist $c', R'>0$ such that 
\begin{equation}
|A_r|\geq c'\,r^{d-s}
\label{euni2.4}
\end{equation}
for all $r\in\langle0,R'\rangle$.
\end{prop}
\proof\
It suffices to prove that there are $b'>0$ and $R'\in\langle0,R/2\rangle$ such that
\begin{equation}
|\Omega_{z, R}\cap B(x\,; r)|\geq b'\,r^d
\label{euni2.41}
\end{equation}
for all $x\in\Gamma_{\!z,R/2}$ and $r\in\langle0,R'\rangle$.
Then  Lemma~\ref{luni2.2}.\ref{luni2.2-2} applied to $A= \Gamma_{\!z,R/2}$ gives
\[
|A_r|=|\{x\in \Omega:d_A(x)< r\}|=| \{x\in\Omega_{z,R}:d_{\Gamma_{\!z,R/2}}(x)<r\}|\geq c'\,r^{d-s}
\]
for $r\in\langle0,R'\rangle$.

The estimate (\ref{euni2.41}) is, however, a consequence of  the  uniformity of $\Omega_{z,R}$ by the following argument.

 Fix $x\in \Gamma_{\!z,R/2}$ and $y\in \Omega$ with $d(x\,;y)=r\in\langle0,R/2\rangle$.
 Therefore $y\in \Omega_{z,R}$.
 Let  $\gamma\colon[0,1]\mapsto \Omega\cup\{x\}$ be a curve joining $x$ and~$y$ which satisfies the local uniformity properties.
 Then $|\gamma(x\,;y)|\leq \sigma\,d(x\,;y)=\sigma \,r$.
 Thus  $d(t\,;x)\leq \sigma \,r$ for each $t\in\gamma([0,1])$.
 But one also has $d_{\Gamma}(t)\geq \sigma^{-1}(d(t\,;x)\wedge d(t\,;y))$.
 In particular if $t$ is the midpoint of $\gamma$ then $2\,d(t\,;x)\geq r$ and
 $d_{\Gamma}(t)\geq (2\sigma)^{-1}r$.
 Consequently, $B(t\,;(2\sigma)^{-1}r)\subseteq B(x\,; (\sigma+(2\sigma)^{-1})\,r)\cap\Omega$.
 Then replacing $r$ by $(\sigma+(2\sigma)^{-1})^{-1}\,r$ one deduces that 
 $B(t\,;\rho\,r)\subseteq \Omega\cap B(x\,; r)$ for all $r\in\langle0,R'\rangle$ with $\rho=(1+2\sigma^2)^{-1}$
 and $R'=2\sigma\rho\,R$.
 Then (\ref{euni2.41}) follows immediately with $b'=|B(0\,;1)|\,\rho^d$.
 \hfill$\Box$

\bigskip

The bounds   on $|A_r|$ for the subsets $A=\Gamma_{\!x,R}=\Gamma \cap B(x\,;R)$ of the boundary $\Gamma$ of $\Omega$ are fundamental for the proof of Theorem~\ref{tuni1.1}.
The estimates are also related to the existence of the Minkowski dimension $d_M(A)$ of $A$.
There are a number of possible definitions of the Minkowski dimension  but the appropriate definition in the current context would be
 \[
 d_M(A)=d-\lim_{r\to0}\log|A_r|/\log r
 \]
 whenever the limit exists.
It follows, however,  from the Ahlfors $s$-regularity of $\Gamma$ that  the limit exists and $d_M(A)=d_H(A)=s$.
(For a fuller  discussion of Ahlfors regularity property and the equality  of various possible dimensions   see  \cite{Leh3}, Lemma~2.1, and \cite{LeT}, Theorems~4.1 and 4.2).

\smallskip

Next  we note that the Ahlfors $s$-regularity property (\ref{euni1.3})  implies local equivalence of Hausdorff measure and Hausdorff content.
This is directly related to the observation that the regularity property implies local equivalence of the measure $\mu$ and the Hausdorff measure $\ch^s$.
In fact the lower bound  in the Ahlfors property (\ref{euni1.3}) implies that if $A=\Gamma\cap B(x_0\,;R_0)$ with $x_0\in\Gamma$ and $R_0>0$ then 
there is an $a>0$ such that $\ch^s(E)\leq a\,\mu(E)$ for all Borel subsets $E\subseteq A$.
Conversely, the upper bound of (\ref{euni1.3})  implies that there is a $b>0$ such that $\mu(E)\leq b\,\ch^s(E)$ for all $E\subseteq A$.
Hence $\mu$ and $\ch^s$ are equivalent on $A$ and $s=d_H(A)$.
(See, for example, \cite{Hei}, Section~8.7.)

Now the Hausdorff measure $\ch^s(E)$ of each Borel subset $E$ of $\Ri^d$ is defined by
 \[
 \ch^s(E)=\lim_{t\to0}\ch^s_t(E)
 \]
 where 
 \[
\ch^s_t(E)=\inf\Big\{\sum_{j=1}^\infty\diam(U_j)^s : \diam(U_j)< t\,;\;E\subseteq \bigcup_jU_j\Big\}\;,
\]
for all $s,t>0$.
But the  Hausdorff content  of $E$ is defined as $\ch^s_\infty(E)$, i.e.\ there is no restriction on the diameters of the
sets in the cover.
Moreover, in these definitions it suffices to consider covers of $E$ by balls $B_j=B(x_j\,;r_j)$ with $x_j\in E$.

\begin{lemma}\label{luni2.40}
If  $\,\Gamma$ satisfies the Ahlfors $s$-regularity property $(\ref{euni1.3})$ and  $A=\Gamma\cap B(x_0\,;R_0)$ with $x_0\in\Gamma$ and $R_0>0$
then there is a $c>0$ such that 
\[
\ch^s_\infty(E)\leq \ch^s(E)\leq c\,\ch^s_\infty(E)
\]
for all Borel subsets $E\subseteq A$.
\end{lemma}
\proof\
The lower bound follows directly from the definition of $\ch^s$ and $\ch^s_\infty$.
To establish the upper bound let $\{B_j\}_{1\leq j\leq N}$ be a covering of $E$ by balls $B_j=B(x_j\,;r_j)$ with $x_j\in E$ and $r_j>0$.
Then 
\[
\ch^s(E)\leq a\, \mu(E)\leq a\sum_{j=1}^N\mu(B_j)\leq a\,c\sum_{j=1}^N r_j^s
\]
and taking the infimum over the possible covers one deduces that $\ch^s(E)\leq a\,c\,\ch^s_\infty(E)$.
\hfill$\Box$

\bigskip

Finally we derive an estimate which  is relevant to the derivation of a local version of the weighted Hardy inequality.
This will
be of importance in the sequel.

\begin{prop}\label{puni2.4}
Assume the boundary $\Gamma$ of the domain $\Omega$ satisfies the  Ahlfors $s$-regularity condition $(\ref{euni1.3})$.
Fix $z\in\Gamma$ and $r>0$.
Then there is an $a>0$ such that 
\begin{equation}
\ch^s_\infty(\Gamma\cap B(x\,;2d_\Gamma(x)))\geq a\, d_\Gamma(x)^s
\label{euni2.30}
\end{equation}
for all $x\in\Omega_{z,r}$.
\end{prop}
\proof\
First  if  $y\in \Gamma$ with $d_\Gamma(x)=d(x\,;y)$ then 
$B(y\,;d_\Gamma(x))\subset B(x\,;2d_\Gamma(x))$.
Hence $\ch^s_\infty(\Gamma\cap B(x\,;2d_\Gamma(x)))\geq \ch^s_\infty(\Gamma\cap B(y\,;d_\Gamma(x))$.
Secondly, if $x\in\Omega_{z,r}$ then $\Gamma\cap B(x\,;2d_\Gamma(x))\subseteq \Gamma_{\!z,3r}$.
Thirdly, it follows from the foregoing that the Hausdorff content and  the Hausdorff measure are equivalent on $ \Gamma_{\!z,3r}$.
Therefore  to prove (\ref{euni2.30})  it suffices to prove that there is a $b>0$ such
\[
{\cal H}^s (\Gamma\cap B(y\,;r))\geq b\,r^s
\]
for all $y\in \Gamma_{\!z,3r}$.
But since $\ch^s$ is locally equivalent to $\mu$ one has
\[
{\cal H}^s (\Gamma\cap B(y\,;r))\geq a\,\mu(\Gamma\cap B(y\,;r))\geq a\,c\,r^s
\]
for all $y\in \Gamma_{\!z,3r}$ by the regularity assumption.
\hfill$\Box$

\bigskip

\section{Markov uniqueness}\label{S3}

In this section we give the proof of Theorem~\ref{tuni1.1}.
It is in two parts. 
First we prove that the degeneracy bounds $\delta\geq1+(s-(d-1))$  imply Markov uniqueness.
This part of the proof is based on an argument given in \cite{RSi4}.
Secondly, we use local versions of the weighted Hardy inequalities derived in \cite{KL1} \cite{Leh2} to prove that 
Markov uniqueness implies the degeneracy bounds.

The first part of the proof is based on the observation of Proposition~\ref{puni2.1} that Markov uniqueness of $h_0$ is equivalent to the 
property $\capp_{h_0}(\Gamma)=0$.
But it follows from the general monotonicity properties of the capacity that $\capp_{h_0}(\Gamma)=0$ if and only if $\capp_{h_0}(A)=0$ for all
bounded  non-empty  subsets $A$ of $\Gamma$.
The latter property is, however,  a consequence of  the arguments of \cite{RSi4}, Proposition~4.2.

\begin{prop}\label{puni3.1}
Let $A$ be a bounded  non-empty  subset of $\Gamma$.
Assume there are $\delta\geq0$ and $b,R_1>0$ such that 
\begin{equation}
C(x)\leq b\,d_A(x)^{\delta}I
\label{euni3.1}
\end{equation}
for almost all $x\in A_{R_1}$.
Further assume  there is a Borel measure $\mu$ on $A$ and $c,R_2>0$ such that 
\begin{equation}
\mu(A\cap B(x\,;r))\geq c\,r^s
\label{euni3.2}
\end{equation}
for all $x\in A$ and $r\in \langle0,R_2\rangle$.

\smallskip

If  $\delta\geq\delta_c$ where $\delta_c=1+(s-(d-1))$  then $\capp_{h_0}(A)=0$.
\end{prop}
\proof\ 
First  by increasing the value of $b$ and decreasing the value of $c$, if necessary,
one may assume that $R_1=R_2=1$ in Conditions (\ref{euni3.1}) and (\ref{euni3.2}).

Secondly,  it follows from Lemma~\ref{luni2.2}.\ref{luni2.2-1} that there is a $c>0$ such that $|A_r|\leq c\, r^{d-s}$ for all $r\in\langle0,1\rangle$.
This upper bound only uses the regularity bound (\ref{euni3.2}).

Thirdly,  define a sequence of 
 functions $s>0\mapsto  \rho_n(s)\in[0,1]$ by $\rho_n(s)=1$ for all $s\in\langle0, n^{-1}]$, $\rho_n(s)=-\log s/\log n$
for $s\in\langle n^{-1}, 1] $ and $\rho_n(s)=0$ if $s>1$.
Then  set $\eta_{r,n}=\rho_n\circ (r^{-1}d_A)$.
It follows that $0\leq \eta_{r,n}\leq 1$, $\supp\eta_{r,n}\subseteq \overline{A}_r$ and $\eta_{r,n}(x)=1$ if $d_A(x)\leq r/n$.
Therefore to prove that $\capp_{h_0}(A)=0$ it suffices to show that $\inf_{r\in\langle0,1]}\inf_{n\geq1}\|\eta_{r,n}\|_{D(h_N)}=0$.

Since  $0\leq \eta_{r,n}\leq 1$ and $\supp\eta_{r,n}\subseteq \overline{A}_r$ it follows that $\|\eta_{r,n}\|_2^2\leq |A_r|\leq c\,r^{d-s}$ for  $r\leq 1$
where the last estimate uses Lemma~\ref{luni2.2}.\ref{luni2.2-1}.
But  $\eta_{r,n}\in D(h_N)$, by construction, and 
\begin{eqnarray*}
h_N(\eta_{r,n})&=&\int_{A_r}dx\,\Gamma(\eta_{r,n})(x)\\[5pt]
&\leq& b\int_\Omega dx\,d_A(x)^{\delta}\,|(\nabla\eta_{r,n})(x)|^2
\leq b\,(\log n)^{-2}\int_\Omega dx\,\one_{\{x:r n^{-1}\leq d_{A}(x)\leq r\}}\,d_A(x)^{-(2-\delta)}
\end{eqnarray*}
where we have used 
the degeneracy bounds 
(\ref{euni3.1}).
If $\delta\geq 2$ then $h_N(\eta_{r,n})\leq b\,|A_1|\,(\log n)^{-2}$
and $\inf_{n\geq 1}\|\eta_{r,n}\|^2_{D(h_N)}\leq c\,r^{d-s}$.
Since this conclusion  holds for all small $r$, and $d>s$,  one deduces that
$\inf_{r\in\langle0,1]}\inf_{n\geq1}\|\eta_{r,n}\|_{D(h_N)}=0$. 
Hence $\capp_{h_0}(A)=0$.

If, however,  $\delta<2$ then
\begin{eqnarray*}
\int_\Omega dx\,\one_{\{x: rn^{-1}\leq d_{A}(x)\leq r\}}\,d_A(x)^{-(2-\delta)}
&\leq&r^{\delta-2} \int_\Omega dx\,\one_{\{x:n^{-1}\leq r^{-1}d_A(x)\leq 1\}}\,(r^{-1}d_A(x))^{-(2-\delta)}\\[5pt]
&&\hspace{-3cm}=r^{\delta-2}  \int_\Omega dx\,\one_{\{x:n^{-1}\leq r^{-1}d_A(x)\leq 1\}}\,\Big(1+(2-\delta)\int^1_{r^{-1}d_A(x)}dt\,t^{-(3-\delta)}\Big) 
\\[5pt]
&&\hspace{-3cm}\leq r^{\delta-2} \Big(|A_{r}|+(2-\delta)\int^1_{n^{-1}} dt\,t^{-(3-\delta)}|A_{rt}|\Big)\;.\label{euni3.0}
\end{eqnarray*}
But  $|A_r|\leq c\,r^{d-s}$ and $|A_{rt}|\leq c\,(rt)^{d-s}$ by Lemma~\ref{luni2.2}.\ref{luni2.2-1}.
Therefore 
\begin{eqnarray*}
\int_\Omega dx\,\one_{\{x: rn^{-1}\leq d_{A}(x)\leq r\}}\,d_A(x)^{-(2-\delta)}&\leq& c\,r^{\delta-\delta_c}
\Big(1+(2-\delta)\int^1_{n^{-1}} dt\,t^{-1} t^{\delta-\delta_c}\Big)\\[5pt]
&\leq&c\, (1+(2-\delta)\log n)
\end{eqnarray*}
where we have used the assumptions $r\leq 1$ and  $\delta\geq \delta_c$.
It follows by combination of these estimates  that  if $2>\delta\geq \delta_c$ then
\[
\|\eta_{r,n}\|^2_{D(h_N)}\leq c\,r^{d-s}+b\,c\,(2-\delta)\,(\log n)^{-1}+b\,c\,(\log n)^{-2}
\]
for all $r\in\langle0,1]$.
Therefore $\inf_{r\in\langle0,1]}\inf_{n\geq1}\|\eta_{r,n}\|_{D(h_N)}=0$
and  $\capp_{h_0}(A)=0$.\hfill$\Box$

\bigskip

The  assertion in Theorem~\ref{tuni1.1} that the bound $\delta\geq \delta_c$ suffices to establish Markov uniqueness is now a corollary of Propositions~\ref{puni2.1} and \ref{puni3.1}.
First note that $d_A(x)\geq d_{\Gamma}(x)$ for all $x\in\Omega$ so the bounds (\ref{euni3.1}) formulated with $d_A$ follow from the similar bounds
formulated with $d_{\Gamma}$.
Therefore the upper bound of the degeneracy condition (\ref{euni1.30}) is sufficient to deduce from
Proposition~\ref{puni3.1} that $\capp_{h_0}(A)=0$ for all the  subsets $A=\Gamma \cap B(x\,;R)$ with $x\in\Gamma$.
Then $\capp_{h_0}(\Gamma)=0$ by the monotonicity properties of the capacity.
Finally $h_D=h_N$ by Proposition~\ref{puni2.1}.

\smallskip

\begin{remarkn}\label{runi3.1}$\;$ 
Proposition~\ref{puni3.1} can be  strengthened by  comparison of the capacity $\capp_{h_0}$ as a non-additive measure on $\Gamma$ and the 
Hausdorff measure.
The argument adapts well known results for the Laplacian on $\Ri^d$ and the classical capacity (see, for example,
\cite{EvG},  Section~4.7   or  \cite{MZ1}, Section~2.1.7).
This approach was used  in \cite{RSi5} Proposition~4.4.
In particular if $\delta<2$ one obtains a bound $\capp_{h_0}(A)\leq c\,\ch^{d+\delta-2}(A)$.
But the regularity assumption (\ref{euni3.2}) implies that $d_H(A)\leq s$.
Therefore if $\delta>\delta_c$ then $d_H(A)\leq s< d+\delta-2$ and $\ch^{d+\delta-2}(A)=0$.
Hence $\capp_{h_0}(A)=0$.
If $\delta=\delta_c$ the argument is slightly more intricate.
Then the regularity property implies that $\ch^{d+\delta-2}(A)<\infty$ and this suffices to deduce that $\capp_{h_0}(A)=0$
by adapting the reasoning of Section~4.7.2 in \cite{EvG} or Theorem~2.52 in \cite{MZ1}.
\end{remarkn}

Next we turn to the proof of the converse statement in Theorem~\ref{tuni1.1}, the assertion that Markov uniqueness of the form $h_0$ implies that $\delta\geq \delta_c$.
The proof is based on weighted Hardy inequalities which are local versions of the Hardy inequalities given by Theorem~1.4 of \cite{KL1} and
Theorem~1.2 of~\cite{Leh2}. 
In conformity with  these references we state the following propositions for all $p\in\langle1,\infty\rangle$ 
although  in the current context they are only of interest for the case $p=2$.
We begin with the case where the degeneracy parameter $\delta$, i.e.\ the weight exponent,  does not exceed $p-1$. 
Again  $\Omega_{z,r}=\Omega\cap B(z\,;r)$ for $z\in\Gamma$ and $r>0$.

\begin{prop}\label{puni3.2*} Let $\Omega$ be a domain in $\Ri^d$ with boundary $\Gamma$.
Fix $z\in\Gamma$ and $r>0$.
Assume  there exist $a>0$ and $s\in[0,d-1]$ such that 
\begin{equation}
\ch^s_\infty(\Gamma\cap B(x\,;2d_\Gamma(x)))\geq a\, d_\Gamma(x)^s
\label{euni3.3*}
\end{equation}
for all $x\in\Omega_{z,r}$.

Then for each $p\in\langle1,\infty\rangle$ and $\delta < p+s-d$ there exists $b>0$ such that the local weighted Hardy inequality
\begin{equation}
\int_{\Omega_{z,r}}dx\,d_\Gamma(x)^\delta\,|(\nabla\varphi)(x)|^p
\geq b \int_{\Omega_{z,r}}dx\, d_\Gamma(x)^{\delta-p}\,|\varphi(x)|^p
\label{euni3.4}
\end{equation}
is valid for all $\varphi\in C_c^\infty(\Omega_{z,r})$.
\end{prop}
\proof\
The proposition is essentially a corollary of the proof of Theorem~1.2 of~\cite{Leh2}.
Assume first that $\delta < p+s-d$ and $\delta\le 0$.
Then it follows from the assumptions and Theorem~4.2 of~\cite{Leh2} that one has a pointwise version
of the Hardy inequality~(\ref{euni3.4}) for all $x\in\Omega_{z,r}$ and $\varphi\in C_c^\infty(\Omega)$.
Explicitly there are $q\in\langle1,p\rangle$ and $c>0$ such that 
\[
|\varphi(x)|\leq c\,d_\Gamma(x)^{1-\delta/p}\Big(M_{2d_\Gamma(x)}(d_\Gamma(x)^{\;\delta q/p}\,|(\nabla\varphi)(x)|^q\Big)^{1/q}
\]
for all $\varphi\in C_c^\infty(\Omega)$ and $x\in\Omega_{z,r}$.
Here $M_R$ denotes the restricted maximal function defined by
\[
(M_R\psi)(x)=\sup_{0<r<R}|B(x\,;r)|^{-1}\int_{B(x;r)}dy\,|\psi(y)|
\]
where $R$ is allowed  to depend on $x$.
But the maximal function is $L^{p/q}$-bounded.
Therefore there is a $c_{p/q}>0$ such that 
\begin{eqnarray*}
 \int_{\Omega_{z,r}}dx\, d_\Gamma(x)^{\delta-p}\,|\varphi(x)|^p&\leq&
  c\int_{\Omega_{z,r}}dx\,\Big(M_{2d_\Gamma(x)}(d_\Gamma(x)^{\;\delta q/p}\,|(\nabla\varphi)(x)|^q\Big)^{p/q}\\[5pt]
  &\leq&c\, c_{p/q}\int_{\Omega_{z,r}}dx\,\Big(d_\Gamma(x)^{\delta q/p}\,|(\nabla\varphi)(x)|^q\Big)^{p/q}\\[5pt] 
  &=&b\int_{\Omega_{z,r}}dx\,d_\Gamma(x)^{\delta }\,|(\nabla\varphi)(x)|^p
  \end{eqnarray*}
  for all $\varphi\in C_c^\infty(\Omega_{z,r})$ where $b=c\, c_{p/q}$.
This proves the claim provided that $\delta\le 0$.

For $0<\delta< p+s-d \le p+(d-1)-d = p-1$ we use the following fact which is essentially 
from Lemma~2.1 of \cite{Leh4}: If the Hardy inequality~(\ref{euni3.4}) is valid 
for all $\varphi\in C_c^\infty(\Omega_{z,r})$ with parameters $\tilde p$ and $\tilde \delta$, and if
$\alpha>0$, then a corresponding Hardy inequality also holds 
for all $\varphi\in C_c^\infty(\Omega_{z,r})$
 with parameters
$\tilde p+\alpha$ and $\tilde \delta+\alpha$. 
This result is formulated in~\cite{Leh4} only for
functions in $C_c^\infty(\Omega)$, but the same proof applies verbatim 
for functions in $C_c^\infty(\Omega_{z,r})$. Now, since $p-\delta>1$ 
and $0 < (p-\delta) +s-d$, we have by the first part of the proof that
the Hardy inequality~(\ref{euni3.4}) is valid 
for all $\varphi\in C_c^\infty(\Omega_{z,r})$ with parameters $\tilde p = p-\delta$ and $\tilde \delta=0$.
Then choosing $\alpha=\delta$ we deduce the Hardy 
inequality~(\ref{euni3.4}) for $p$ and $\delta$ from Lemma~2.1 in \cite {Leh4}. \hfill$\Box$
  
  \bigskip

For weight exponents $\delta\ge p-1$, the thickness condition~(\ref{euni3.3*})
alone is not sufficient for the Hardy inequality~(\ref{euni3.4}), since 
in this case the geometry of the boundary $\Gamma$ also affects the
validity of Hardy inequalities. Here we follow the ideas in~\cite{KL1}
and use a local version of the sufficient condition for
Hardy inequalities 
formulated in terms of the visual boundary near a point $x\in \Omega$.
For a fixed $\tau\geq 1$
this set, which we denote by $\Gamma_{\rm\! \!vis}(x)$, consists of those  $y$  in the boundary $\Gamma$ of $\Omega$ for which there is 
rectifiable curve $\gamma\colon[0,1]\mapsto \Omega\cup\{y\}$ with $\gamma(0)=y$ and $\gamma(1)=x$ satisfying
$d_\Gamma(\gamma(t))\geq \tau^{-1} d(\gamma(t) ; y)$ for all $t\in[0,1]$.
Although it is not generally true that $\Gamma_{\rm \!\!vis}(x)\subset \Omega_{z,r}$ for $x\in\Omega_{z,r}$ it is nevertheless true that there is an $L\geq1$, whose value depends on $\tau$, such that 
$\Gamma_{\rm \!\!vis}(x)\subset \Omega_{z,Lr}$ for all $x\in \Omega_{z,r}$.

Now we have the following localized version of Theorem~1.4 of \cite{KL1}.

\begin{prop}\label{puni3.2} Let $\Omega$ be a domain in $\Ri^d$ with boundary $\Gamma$.
Fix $z\in\Gamma$ and $r>0$.
Assume 
there exist $a>0$ and $s\in[0,d]$ such that 
\begin{equation}
\ch^s_\infty(\Gamma_{\rm\!\!vis}(x))\geq a\,d_\Gamma(x)^s
\label{euni3.3}
\end{equation}
for all $x\in\Omega_{z,r}$.

Then for each $p\in\langle1,\infty\rangle$ and $\delta < p+s-d$ there exists a $b>0$ such that the 
local weighted Hardy inequality~$(\ref{euni3.4})$
is valid for all $\varphi\in C_c^\infty(\Omega_{z,r})$.
\end{prop}
\proof\
The proposition follows from the 
the proof of Theorem~1.4 of \cite{KL1} in
exactly the same way as Proposition~\ref{puni3.2*}
was proven above. Indeed,
from the assumptions and Theorem~5.1 of~\cite{KL1} one  again obtains a pointwise version
of the Hardy inequality for all $x\in\Omega_{z,r}$ and $\varphi\in C_c^\infty(\Omega)$.
Then the Hardy inequality~(\ref{euni3.4}) follows  for all $\varphi\in C_c^\infty(\Omega_{z,r})$ 
from the $L^{p/q}$-boundedness of the maximal function 
in the same way as in the proof of Proposition~\ref{puni3.2*}.\hfill$\Box$

  \bigskip
  
If $p=2$ then the  Ahlfors regularity and local uniformity properties  introduced in Section~\ref{S1} allow
one to deduce the following version of the Hardy inequality.

\begin{cor}\label{cuni3.1}
Fix $z\in\Gamma$ and $r>0$.
Assume $A=\Gamma\cap B(z\,;3r)$ satisfies the Ahlfors $s$-regularity property $(\ref{euni1.3})$ and that $\delta < 2 + s -d = 1 + (s-(d-1))$.
If $\delta\geq 1$   we also assume that $\Omega_{z,r}$ is $\Omega$-uniform.

Then there exists an $a>0$ such that 
\begin{equation}
\int_{\Omega_{z,r}}dx\,d_\Gamma(x)^\delta\,|(\nabla\varphi)(x)|^2
\geq a \int_{\Omega_{z,r}}dx\, d_\Gamma(x)^{\delta-2}\,|\varphi(x)|^2
\label{euni3.5}
\end{equation}
is valid for all $\varphi\in C_c^\infty(\Omega_{z,r})$.
\end{cor}
\proof\
If $s\le d-1$, 
the corollary is an immediate consequence of Proposition~\ref{puni3.2*}.
Notice that in this case necessarily $\delta<1$.

If $s>d-1$ but $\delta<1$, we use the fact that
$s$-regularity implies the thickness condition (\ref{euni3.3*})
for $s$ (cf.\ Proposition~\ref{puni2.4}), and hence
also for the exponent $d-1$.
Then the claim follows again from Proposition~\ref{puni3.2*}
applied with $s=d-1$.

Finally, if $s>d-1$ and $\delta\ge 1$, we use Proposition~\ref{puni3.2}.
Therefore we need to know that the visual boundary condition (\ref{euni3.3})
holds for all $x\in \Omega_{z,r}$.
But the local uniformity condition ensures that $\Gamma\cap B(x\,;2d_\Gamma(x))\subset \Gamma_{\rm\!\!vis}(x)$
(see Proposition~4.3 in \cite{KL1}).
Therefore
\[
\ch^s_\infty(\Gamma_{\rm\!\!vis}(x))\geq \ch^s_\infty(\Gamma\cap B(x\,;2d_\Gamma(x)))\geq a\, d_\Gamma(x)^s
\]
by Proposition~\ref{puni2.4}. Hence the assumption of Proposition~\ref{puni3.2} is valid,
and the Hardy inequality follows.
\hfill$\Box$

\bigskip

Note that  the factor $3$ occurring in  the corollary in the  radius of the enlarged balls has no particular significance.
It could be replaced by any $\lambda>1$.

\begin{prop}\label{puni3.4}
Fix $z\in\Gamma$ and $r>0$.
Assume $A=\Gamma\cap B(z\,;3r)$ satisfies the Ahlfors $s$-regularity property $(\ref{euni1.3})$ with $s\in\langle0,d\rangle$ 
and that there are $b,\delta>0$ such that 
\begin{equation}
C(x)\geq b\,d_\Gamma(x)^\delta I
\label{euni3.6}
\end{equation}
for almost all $x\in \Omega_{z,r}$.
If $\delta\geq 1$   we also assume that $\Omega_{z,r}$ is $\Omega$-uniform.

\smallskip

It follows that if $h_0$ is Markov unique then $\delta\geq 1+(s-(d-1))$.
\end{prop}
\proof\
The proof is in three steps.
First  the lower bound (\ref{euni3.6}) on the coefficients of the form
$h_0$ together with the Hardy inequality (\ref{euni3.5})
 gives the bounds
 \begin{equation}
 h_0(\varphi)\geq b\, c \int_{\Omega_{z,r}}dx\, d_\Gamma(x)^{\delta-2}\,|\varphi(x)|^2
 \label{euni3.7}
 \end{equation}
 for all $C_c^\infty(\Omega_{z,r})$.
Secondly,  the Markov uniqueness assumption allows the extension of  (\ref{euni3.7}) by continuity to a positive $\varphi\in D(h_N)$
with $\supp\varphi\subset \Omega_{z,r}$ and $\varphi=1$ on $\Omega_{z,\rho}$ where $\rho\in\langle0,r\rangle$.
One then immediately deduces that 
\[
\int_{\Omega_{z,\rho}}dx\,  d_\Gamma(x)^{\delta-2}<\infty
\;.
\]
Thirdly, one argues that if $\delta<2+s-d$ then there is a  contradiction.
Hence one must have $\delta\geq 2+s-d$

\smallskip
The details of the last two steps are as follows.
\smallskip

Fix $\varphi\in C_c^\infty(\Ri^d)$ with $0\leq \varphi\leq 1$, $\supp\varphi\subset B(z\,;\rho+\tau)$ and $\varphi=1$ on $B(z\,;\rho)$ where $\rho,\tau>0$
and $\rho+\tau<r$.
Then the restriction $\varphi|_\Omega$  of $\varphi$ to $\Omega$, which we identify with $\varphi$, is in $D(h_N)$.
But since $h_N=h_D$ there is a sequence $\varphi_n\in C_c^\infty(\Omega)$ which converges to $\varphi$ in the $D(h_N)$-graph
norm.
Now fix $\chi\in C_c^\infty(B(z\,;r))$ such that $0\leq \chi\leq 1$ and $\chi=1$ on $B(z\,;\rho+\tau)$.
Again  identifying $\chi$ with $\chi|_\Omega$ one has $\chi\,\varphi=\varphi$.
But, by Leibniz' rule, 
\begin{eqnarray*}
\|\chi\,\varphi_n-\varphi\|_{D(h_N)}^2&=&\|\chi\,(\varphi_n-\varphi)\|_{D(h_N)}^2\\[5pt]
&\leq& 2\,\|\chi\|_\infty \,h_N(\varphi_n-\varphi) +(2\,\|\Gamma(\chi)\|_\infty+\|\chi\|_\infty^2)\,\|\varphi_n-\varphi\|_2^2
\;.
\end{eqnarray*}
Therefore $\|\chi\,\varphi_n-\varphi\|_{D(h_N)}\to0$ as $n\to\infty$.
Hence $\varphi$ is also approximated in the $D(h_N)$-graph norm by 
the $\chi\,\varphi_n\in C_c^\infty(\Omega_{z,r})$.
Thus 
\[
\infty>h_N(\varphi)\geq b\, c \int_{\Omega_{z,r}}dx\, d_\Gamma(x)^{\delta-2}\,|\varphi(x)|^2
\geq b\,c\int_{\Omega_{z,\rho}}dx\, d_\Gamma(x)^{\delta-2}
\]
by  (\ref{euni3.7}) and a continuity argument.

Finally let $\Omega_{z,\rho:\tau}=\Omega_{z,\rho}\cap \Gamma_\tau$
where $\Gamma_\tau$ is again the $\tau$-neighbourhood of $\Gamma$.
If $\delta< 2+s-d$
then $\delta<2$ and there is an $a>0$ such that
\[
\int_{\Omega_{z,\rho}}dx\, d_\Gamma(x)^{\delta-2}\geq
\int_{\Omega_{z,\rho:\tau}}dx\, d_\Gamma(x)^{\delta-2}\geq  \tau^{-(2-\delta)}\,|\Omega_{z,\rho:\tau}|\geq a\,\tau^{-(2+s-d-\delta)}
\]
 for $\rho$ fixed and uniformly for all $\tau>0$ with  $\rho+\tau<r$.
 The second step uses the lower bound of  Proposition~\ref{puni2.2} if $s\in\langle0,d-1\rangle$ or the lower bound of  Proposition~\ref{puni2.3} if  $s\in[ d-1,d\rangle$.
 Since  the lower bound on the integral diverges as $\tau\to0$ one obtains a contradiction.
 Therefore  $\delta\geq 2+s-d=1+(s-(d-1))=\delta_c$.
 \hfill$\Box$

\bigskip

\noindent{\bf Proof of Theorem~\ref{tuni1.1}}$\;$ The theorem follows by combination of Propositions~\ref{puni3.1} and \ref{puni3.4}.
The degeneracy bounds (\ref{euni1.30}) include the upper (\ref{euni3.1}) and lower (\ref{euni3.6}) bounds on the coefficients.
Hence both propositions are applicable.
\hfill$\Box$

\bigskip

The foregoing proofs are essentially local and depend only on  properties in the neighbourhood of the 
boundary.

The proof that the degeneracy bound $\delta\geq \delta_c$ implies Markov uniqueness is  a consequence of the local capacity estimates $\capp_{h_0}(A)=0$ for a suitable  family of bounded subsets $A$ of $\Gamma$.
It  depends on the assumption that $\delta$ is constant, i.e.\ the  value of $\delta$  is independent of  the choice of $A$.
The neighbourhood of $\Gamma$ enters the proof through the  estimates $|A_r|\leq b\,r^{d-s}$ on the $r$-neighbourhood $A_r$ of $A$.
These estimates are independent of the local $\Omega$-uniformity (see Lemma~\ref{luni2.2}.\ref{luni2.2-1}).

The proof that Markov uniqueness implies $\delta\geq \delta_c$ is strictly local; it only requires estimates on a neighbourhood $\Omega_{z,r}$ for one $z\in \Gamma$ and all small $r>0$.
The global nature of the conclusion follows because $\delta$ and $s$ are constant.
If $s\geq d-1$ then the $\Omega$-uniformity of $\Omega_{z,r}$ enters the proof   in  two distinct ways.
First it is used to establish the local Hardy inequality through local identification of the visual boundary.
Secondly, it is used to derive the local lower bounds $|A_r|\geq c\, r^{d-s}$ with $A=\Gamma\cap B(z\,;\rho)$.
Note that the second step in the proof of Proposition~\ref{puni3.4} can be reformulated in terms of the capacity.
Markov uniqueness is equivalent to $\capp_{h_0}(\Gamma)=0$ and this requires that $\capp_{h_0}(\Gamma_{\!z,r})=0$.
But a slight variation of the argument with the local Hardy inequality establishes that the latter capacity condition 
is incompatible with the degeneracy condition $\delta<2+s-d$.

It might appear surprising that in the latter proof one only needs estimates near one point $z\in\Gamma$
and in the case $s\geq d-1$ this has to be a point of local uniformity.
This can, however, be  understood by noting that Markov uniqueness  implies that the 
boundary $\Gamma$ is inaccessible to the diffusion and this means that all parts of the boundary must be inaccessible.
But the  points $z\in\Gamma$ at which there is an $\Omega$-uniform neighbourhood $\Omega_{z,r}$ are  potentially the most accessible.
Therefore if the degeneracy of the coefficients is sufficient to ensure that the corresponding sections $\Gamma_{\!z,r}$
of the boundary are not accessible to the diffusion then the rest of the boundary is automatically inacessible.
This is the essence of the proof.
The condition $s\geq d-1$ is significant since it includes the case that the topological dimension of the boundary is $d-1$.
In the latter situation the boundary separates $\Omega$ from its complement and provides a substantial barrier to the diffusion.
In the low dimensional case, $s<d-1$, the boundary is relatively negligible and uniformity is unnecessary.

 \smallskip

These observations  immediately lead to a more general result.

 \begin{cor}\label{cuni4.1}
Let  $\Omega$ satisfy the assumptions of Theorem~$\ref{tuni1.1}$ and let $\Omega'=\Omega\backslash\Gamma'$
where $\Gamma'\subset \Omega$ is a closed  Ahlfors $s'$-regular set  with  $s'\in\langle0,s]$.
Further assume   the coefficients of the form $h_0$ satisfy the degeneracy condition $(\ref{euni1.30})$ on the boundary $\partial\Omega'=\Gamma\cup \Gamma'$ of
$\Omega'$ with 
$\delta\geq 0$.
 
\smallskip

 Then the  form $h_0$ is Markov unique if and only if 
  $\delta\geq 1+(s-(d-1))$.
 \end{cor}
 \proof\
 First if $s'=s$ then  the assumptions of Theorem~$\ref{tuni1.1}$ are satisfied with $\Omega$ replaced
 $\Omega'$  and $\Gamma$ replaced  by $\Gamma\cup\Gamma'$.
 Therefore there is nothing to prove.
Secondly, if $s'<s$ then  $\delta\geq 1+(s-(d-1))$ implies   $\delta> 1+(s' -(d-1))$.
Hence one deduces that $\capp_{h_0}(\Gamma)=0=\capp_{h_0}(\Gamma')$ from Proposition~\ref{puni3.1}.
Therefore $\capp_{h_0}(\partial\Omega')=\capp_{h_0}(\Gamma)+\capp_{h_0}(\Gamma')=0$ and $h_D=h_N$
 by Proposition~\ref{puni2.1}.
 
Thirdly, if $h_D=h_N$ then it follows from Proposition~\ref{puni3.4} that $\delta\geq 1+(s-(d-1))$.
 Note that as $\Gamma$ and $\Gamma'$ are disjoint the local arguments of the proposition 
 still apply to $\Gamma$ if $r$ is sufficiently small.\hfill$\Box$

 \bigskip
  
One can also prove an analogue of   Corollary~\ref{cuni4.1} in  which a countable subset of $\Omega$ is excised.
This is essentially an $s'=0$ version of the foregoing.

 \begin{cor}\label{cuni4.10}
 Let  $\Omega$ satisfy the assumptions of Theorem~$\ref{tuni1.1}$ and let $\Omega'=\Omega\backslash\Gamma'$
where $\Gamma'$ is a countable subset of $\Omega$.
Set $h_0'=h_0|_{C_c^\infty(\Omega')}$.
\smallskip

Then the  form $h_0'$ is Markov unique if and only if $d\geq 2$ and  $\delta\geq 1+(s-(d-1))$.
\end{cor}

Note that in this corollary the coefficients of $h_0$ are assumed to satisfy the degeneracy condition (\ref{euni1.30}) on the boundary $\Gamma$ of $\Omega$
but  not  on the excised set $\Gamma'$.

\smallskip

\noindent{\bf Proof of Corollary~\ref{cuni4.10}}
First assume the regularity condition and the degeneracy condition    hold on $\Gamma$. 
Then it follows as before  that 
$\capp_{h_0}(\Gamma)=0$. 
Since the estimates are localized on the boundary $\Gamma$ it also  follows that $\capp_{h_0'}(\Gamma)=0$. 
Next we argue that if $d\geq 2$ then $\capp_{h_0'}(\Gamma')=0$.

Let $A=\{x_k\}$ with $x_k\in \Gamma'$. 
Then $|A_r|\leq c\,r^d$ for all small $r$.
Consequently,  if $d\geq 2$  one concludes that $\capp_{h_0'}(A)=0$ by the estimates in the proof of Proposition~\ref{puni3.1}.
(Effectively $s=0=\delta$  on $\Gamma'$   and then $d\geq 2$ is the special case of the condition $\delta\geq 1+(s-(d-1))$ used in the calculation.)
Since this argument applies for all $x_k\in \Gamma'$ it follows that $\capp_{h_0'}( \Gamma')=0$ whenever $d\geq 2$. 
Therefore $\capp_{h_0'}(\partial\Omega')=\capp_{h_0'}(\Gamma)+\capp_{h_0'}( \Gamma')=0$ and $h_0'$ is Markov unique by Proposition~\ref{puni2.1}.

Conversely if  $h_0'$ is Markov unique one must have $\capp_{h_0'}(\partial\Omega')=0$.
But since $h_0$ is an extension of $h_0'$ this implies that  $\capp_{h_0}(\partial\Omega')=0$.
In particular  $\capp_{h_0}(\partial\Omega)=0$ and one deduces that  $\delta\geq 1+(s-(d-1))$ by the earlier arguments with the local Hardy inequality. 
Next the Markov uniqueness also implies that $\capp_{h_0'}(\Gamma')=0$.
Therefore  $\capp_{h_0'}(\{x_k\})=0$ for all $x_k\in\Gamma'$.
Now let $B_k=B(x_k\,;r)$ where $r$ is sufficiently small that $B_k\subset \Omega'$ and set $B_k'=B_k\backslash\{x_k\}$.
Since $C(x)>0$ for almost all $x\in\Omega$ it follows that there is a $c_k>0$ such that $h_0(\varphi)\geq c_k\,\|\nabla\varphi\|_2^2$
for all $\varphi\in C_c^\infty(B_k)$.
But there is an $a_k>0$ such that the Hardy inequality 
\[
\|\nabla\varphi\|_2^2\geq a_k\int_{B'_k}dx\,|x-x_k|^{-2}|\varphi(x)|^2
\]
holds for all $\varphi\in C_c^\infty(B_k')$  if $2-d>0$ (see \cite{KL1} Example 7.2).
Therefore 
\[
h_0'(\varphi)\geq c_k\,\|\nabla\varphi\|_2^2\geq  c_k \,a_k\int_{B'_k}dx\,|x-x_k|^{-2}|\varphi(x)|^2
\]
for all $\varphi\in C_c^\infty(B_k')$.
Then as $h_0'$ is Markov unique these bounds extend to all $\varphi\in C_c^\infty(B_k)$.
Now choosing $\varphi$ such that $\varphi=1$ on $B(x_k\,;r/2)$ one has 
\[
\infty> \int_{B(x_k;r/2)}dx\,|x-x_k|^{-2}=\omega_1\int^{r/2}_0ds\,s^{d-3}
\]
with $\omega_1=|B(0\,;1)|$ which contradicts the assumption that $d<2$.
Therefore one must have $d\geq 2$.
\hfill$\Box$

\bigskip

Finally we return to  Corollary~\ref{cuni1.1}.
This  is  a consequence of Theorem~\ref{tuni1.1} and the following lemma. 
Recall that $\Gamma\subset\Ri^d$ is uniformly disconnected, if 
there exists a constant $C\ge 1$ such that for every
$z\in\Gamma$ and all $r>0$ one can find a closed set $A\subset \Gamma$ such that
$\Gamma\cap B(z,r/C)\subset A \subset B(z,r)$ and $\dist(A\,;\Gamma\backslash A)\ge r/C$. 
(For further details and alternative equivalent definitions see \cite{Hei}, Section~14.24.)

\begin{lemma}\label{luni3.1}
Let $\Omega\subseteq \Ri^d$ be a uniform domain and  $\Gamma$   a closed uniformly disconnected subset of $\Omega$.
Then the complement $\Omega\backslash\Gamma$ is also a uniform domain.
\end{lemma}
\proof\ 
MacManus, \cite{McM} page~275, observed that the uniformity of $\Omega'=\Omega\backslash\Gamma$ can proved using the general compactness results of 
V\"ais\"al\"a \cite{Vai} 
but for the sake of completeness we give a direct construction.

Let $x,y\in\Omega'$ and let $\gamma\colon[0,1]\mapsto \Omega$ with $\gamma(0)=x$, $\gamma(1)=y$
be the curve satisfying the uniformity condition with respect to $\Omega$ with the constant $\sigma\ge 1$.
We first modify $\gamma$ into a continuum $E$ which satisfies the conditions for `distance cigars'
in Section~2.4 of~\cite{Vai} with respect to $\Omega'$, with a constant independent of $x$ and $y$.
More precisely, we  show that if we choose $\tilde\sigma=16\,C\sigma>0$, where $C\ge 1$ is
the constant from the  definition of uniform disconnectedness, then
$\diam(E)\le \tilde\sigma d(x\,;y)$ and
$d_{\partial\Omega'}(w)\geq \tilde\sigma^{-1}\lambda(w)$ for all $w\in E$,
where $\lambda(w)=d(x\,;w)\wedge d(w\,;y)$.

If $d_\Gamma(w)\geq \tilde\sigma^{-1}\lambda(w)$ for all 
$w\in\gamma([0,1])$, then we can take $E=\gamma$ and the claim follows. 
Otherwise set $t_0'=0$ and let $t_1$ be the
smallest of the numbers $t\in [t_0',1]$ for which 
$d_\Gamma(\gamma(t)) = \tilde\sigma^{-1}\lambda(\gamma(t))$.
Then set $w_1=\gamma(t_1)$. 
If $d(w_1\,; y) \le d(w_1\,; x)$, 
we can end this part of the construction,
and move to a corresponding construction starting from $y$.
Otherwise, take $z_1\in\Gamma$ with $d_\Gamma(w_1)=d(w_1\,;z_1)$, and let $A_1\subset \Gamma$ be the closed set
given by the definition of uniform disconnectedness for $z_1$ and 
$r_1=4\,C \,d_\Gamma(w_1)=4\,C\,\tilde\sigma^{-1}\lambda(w_1)$. 
Furthermore, let $E_1$ be a connected component of the set
$\{w\in\Ri^d : d_{A_1}(w) = r_1/2C\}$ which intersects $\gamma([0,t_1])$
and separates $x$ and $w_1$;
such a component exists since if $d_{A_1}(w) = r_1/2C$, then
\[
d(x\,; w)\ge d(x\,; w_1) - d(w_1\,; w) \ge \lambda(w_1) - 2r_1 \ge 4r_1 - 2r_1 > r_1/2C
\]
and $d_{A_1}(\gamma(t_1)) = d_{\Gamma}(w_1) < r_1/2C$.

Since $d_{\Gamma}(w) = r_1/2C$ for all $w\in E_1$, we obtain the following estimates: 
$d(w\,;w_1) \le 2r_1$, 
\[
\lambda(w)\le 2r_1 + \lambda(w_1) \le \frac{\tilde\sigma}{4C}r_1 + \frac{\tilde\sigma}{4C}r_1 \le 
\tilde\sigma d_\Gamma(w),
\]
(in the second inequality we used $\tilde\sigma=16\,C\,\sigma \ge 8\,C$), and
\[
d_{\partial\Omega}(w) \ge d_{\partial\Omega}(w_1) - d(w\,;w_1) \ge 
\sigma^{-1} \lambda(w_1) - 2r_1 \ge 2 r_1 \ge \tilde\sigma^{-1} \lambda(w), 
\]
where the penultimate estimate holds since $\sigma^{-1} \lambda(w_1) = 4r_1$ by the choices of
$r_1$ and $\tilde\sigma$. By the above estimates we conclude that
$d_{\partial\Omega'}(w)\ge \tilde\sigma^{-1} \lambda(w)$
for all $w\in E_1$.

Next, let $t_1'$ be the largest of the numbers $t\in [t_1,1]$ for which 
$\gamma(t)\in E_1$. Then $d(\gamma(t_1') \,; w_1) \ge r_1 / 4C =\tilde\sigma^{-1}\lambda(w_1)$.
If $d(\gamma(t_1')\,; y) \le d(\gamma(t_1')\,; x)$, we 
set $t_2=t_1'$ and finish 
this part of the construction. Otherwise we continue inductively and let
$t_2$ be the
smallest of the numbers $t\in [t_1',1]$ for which 
$d_\Gamma(\gamma(t))= \tilde\sigma^{-1}\lambda(\gamma(t))$, and 
denote $w_2=\gamma(t_2)$. If $d(w_2\,; y) \le d(w_2\,; x)$ 
or such a $t_2$ does not exist, we can finish
the construction.
Otherwise, take $z_2\in\Gamma$ with $d_\Gamma(w_2)=d(w_2\,;z_2)$ and let $A_2\subset \Gamma$ be the closed set
given by the definition of uniform disconnectedness for $z_2$ and 
$r_2=4\,C \,d_\Gamma(w_2)$.
Let $E_2$ be a connected component of the set
$\{w\in\Ri^d : d_{A_2}(w) = r_2/2C\}$, such that $E_2$ intersects the set
\[
\gamma([0,t_1]) \cup E_1 \cup \gamma([t_1',t_2])
 \]
and separates $x$ and $w_2$, 
and let $t_2'$ be the largest of the numbers $t\in [t_2,1]$ for which 
$\gamma(t)\in E_2$. As above, we see that $d_{\partial\Omega'}(w)\ge \tilde\sigma^{-1} \lambda(w)$
for all $w\in E_2$ and that $d(\gamma(t_2') \,; w_2) \ge \tilde\sigma^{-1}\lambda(w_2)$.

Continuing this way, we at some point reach $t_n \in  [0,1]$ such that
$d(\gamma(t_n)\,; y) \le d(\gamma(t_n)\,; x)$, since in each step we move from $w_k$
to a point $\gamma(t_k')$ whose distance to $w_k$ is bounded from below by 
$\tilde\sigma^{-1}\lambda(w_k)$.
Now define 
\[
E_x=\bigcup_{k=1}^n  E_{k-1} \cup \gamma([t_{k-1}',t_k])
\]
with  $E_0=\emptyset$.
Then $E_x$ is a continuum joining $x$ to the point $w_n=\gamma(t_n)$.

After this, we make the corresponding construction starting from $y$, i.e.\  from $t=1$, and make $t$ smaller at each
step until we reach some $\tilde t_m \in [0,1]$ such that
$\tilde t_m \le t_n$, where $t_n$ is given above in the construction of $E_x$. We  
obtain the corresponding continuum $E_y$ joining
$y$ to the point $\gamma(\tilde t_m)$. Now $E=E_x\cup E_y$ is a continuum joining $x$ to $y$,
such that $d_{\partial\Omega'}(w)\ge \tilde\sigma^{-1} \lambda(w)$ for all $w\in E$. 
Moreover, it is easy to show that for every $w\in E$, the distance from $w$ to $\gamma$ 
is at most $\diam(\gamma)$, and hence the diameter of $E$ is bounded from above
by $3\diam(\gamma)\le 3\sigma d(x\,;y) \le \tilde\sigma d(x\,;y)$.
Finally, by~\cite{Vai}, Lemma~2.10, we can replace $E$ by a curve 
$\tilde\gamma\colon[0,1]\mapsto \Omega$ with $\tilde\gamma(0)=x$, $\tilde\gamma(1)=y$, 
which satisfies the conditions required in the definition of uniformity for $\Omega'$ 
with the constant $2\,\tilde\sigma$,
and thus $\Omega'$ is indeed uniform.
\hfill$\Box$

\bigskip

\noindent{\bf Proof of Corollary~\ref{cuni1.1}.}$\;$
Assume $\Gamma$ is a uniformly disconnected set and $\Omega=\Ri^d\backslash \Gamma$.
Then $\Omega$ is a uniform set by Lemma~\ref{luni3.1}. 
But $\Gamma$ is the boundary of $\Omega$.
Therefore it  follows from the uniformity of $\Omega$ that if $z\in \Gamma$ and $R>0$ is sufficiently small then
 $\Omega_{z,R}$  is $\Omega$-uniform.
Thus if $\Gamma$ satisfies the Ahlfors $s$-regularity property (\ref{euni1.3}) with $s\in\langle0,d\rangle$ then
Corollary~\ref{cuni1.1} follows directly from Theorem~\ref{tuni1.1}.\hfill$\Box$

\section{General boundaries}\label{S4}

In this section we discuss extensions of Theorem~\ref{tuni1.1} in which the assumptions on the boundary $\Gamma$ are weakened.
The advantage of Theorem~\ref{tuni1.1} is that it covers domains with boundaries of all possible dimensions but the disadvantage is that the Ahlfors
regularity property ensures that the dimension does not vary on the boundary.
This restriction can, however, be relaxed  and the 
 conclusion of Theorem~\ref{tuni1.1} can be extended to domains whose boundaries have various components and faces  with different  regularity properties.
 
 The simplest situation occurs  if the boundary is the union of separated components.

\begin{thm}\label{tuni4.0}
 Let $\Omega$ be a  domain in $ \Ri^d $ whose   boundary $\Gamma$ is the union of a  family $\{\Gamma_{\!\alpha}\}$ of  closed subsets indexed by a possibly uncountable set $M$ where $d(\Gamma_{\!\alpha}\,;\Gamma_{\!\beta})\geq d_0>0$ for all $\alpha, \beta\in M$ with $\alpha\neq \beta$.
Assume that $\Gamma_{\!\alpha}$  is Ahlfors $s_\alpha$-regular with $s_\alpha\in \langle0, d\rangle$.
Further if  $s_\alpha\in[ d-1,d\rangle$  assume there are $z_\alpha\in\Gamma_\alpha$ and $R_\alpha>0$ such that $\Omega_{z_\alpha, R_\alpha}$  is 
$\Omega$-uniform.
Finally assume  the coefficients of the form $h_0$ are $\delta_\alpha$-degenerate on $\Gamma_{\!\alpha}$ 
 for $\delta_\alpha\geq 0$.

\smallskip

 Then the  form $h$ is Markov unique, i.e.\ $h_D=h_N$, if and only if 
  $\delta_\alpha\geq \delta_{c,\alpha}$ for each $\alpha$ where $\delta_{c,\alpha}=1+(s_\alpha-(d-1))$.
 \end{thm}
 \proof\
Note that $d(\Gamma_{\!\alpha}\,;\Gamma_{\!\beta})$ denotes the Euclidean distance between the two boundary components $\Gamma_{\!\alpha}$ and $\Gamma_{\!\beta}$.
Moreover the condition that the coefficients are `$\delta_\alpha$-degenerate on $\Gamma_{\!\alpha}$' 
is understood to mean that for each bounded  non-empty subset $A\subset \Gamma_{\!\alpha}$
there are $a_\alpha,b_\alpha,r_\alpha>0$ such that 
\begin{equation}
a_\alpha\,d_{\Gamma_{\!\alpha}}(x)^{\delta_\alpha}I\leq C(x)\leq b_\alpha\,d_{\Gamma_{\!\alpha}}(x)^{\delta_\alpha}I
\label{euni4.30}
\end{equation}
for  almost all $x\in A_{r_\alpha}$. 
This family of conditions is compatible because of the separation property of distinct components $\Gamma_\alpha$.
Then the proof follows by applying the arguments which established Theorem~\ref{tuni1.1} to each component $\Gamma_{\!\alpha}$.

The proof that $\delta\geq \delta_{c,\alpha}$  implies  $\capp_{h_0}(\Gamma_{\!\alpha})=0$ again follows from applying Proposition~\ref{puni3.1} to an increasing family of 
bounded sets $A\subset \Gamma_{\!\alpha}$. 
The separation assumption ensures that this procedure does not present any additional problems.
Once one has $\capp_{h_0}(\Gamma_{\!\alpha})=0$ for each~$\alpha$ then Markov uniqueness  follows from Proposition~\ref{puni2.1} and the additivity property
$\capp_{h_0}(\Gamma)=\sum_{\alpha\in M}\capp_{h_0}(\Gamma_\alpha)$.

The converse proof that Markov uniqueness implies $\delta_\alpha\geq \delta_{c,\alpha}$  is also evident as it only involves the Hardy inequality 
(\ref{euni3.5}) on the sets $\Omega_{z_\alpha,R_\alpha}=\Omega\cap B(z_\alpha\,;R_\alpha)$ for one $z_\alpha\in \Gamma_{\!\alpha}$ and  one small $R_\alpha>0$.
These inequalities follow, however, from Corollary~\ref{cuni3.1} for each $\alpha$.
We omit further details.
\hfill$\Box$

\bigskip

Theorem~\ref{tuni1.1} also  extends to domains whose boundaries have a finite number of regular faces but with different indices of regularity $s$. 
In the formulation of this extension we endow $\Gamma$ with the relative (Euclidean) topology.

The next result is  divided into two  statements analogous  to Propositions~\ref{puni3.1} and \ref{puni3.4}.

 \begin{thm}\label{tuni4.1}
  Let $\Omega$ be a   domain  with boundary $\Gamma=\bigcup_{j=1}^m \overline{F_j}$ where
the $F_j$ are a finite family of pairwise disjoint Ahlfors $s_j$-regular subsets, with $s_j\in\langle0,d\rangle$, which are  open in the relative topology.
Set $\delta_{c,j}=1+(s_j-(d-1))$.

\begin{tabel}
\item\label{tuni4.1-1}
If  for each $j$ there is  a $\delta_j\geq \delta_{c,j}$  and   $b_j,r_j>0$ such that  for each bounded set $A\subset F_j$ one has   $0\leq C(x)\leq b_j\,d_\Gamma(x)^{\delta_j}I$
for almost all $x\in\Omega$ with $d_A(x)=d_\Gamma(x)<r_j$ then $h_D=h_N$.
 \item\label{tuni4.1-2}
 Further assume  that if   $s_j\geq d-1$  then there is a $z_j\in F_j$ and an $R_j>0$ such that $\Omega_{z_j,R_j}$ is $\Omega$-uniform.
If $h_D=h_N$ and for each $j$ there are $\delta_j\geq 0$ and $b_j,r_j >0$ such that for each bounded set
$A\subset F_j$ one has  $ C(x)\geq b_j\,d_\Gamma(x)^{\delta_j}I$
for almost all  $x\in\Omega$ with $d_A(x)=d_\Gamma(x)<r_j$ then  $\delta_j\geq \delta_{c,j}$.
 \end{tabel}
\end{thm}
\proof\
Although the faces $F_j$ are assumed to be disjoint the relative closures $\overline F_j$ can intersect in lower dimensional `edges' $\overline F_i\cap \overline F_j$,
$\overline F_i\cap \overline F_j\cap\overline F_k$, etc.
This  creates a  new problem in the estimation of the capacity of the various faces.
This is the reason for  considering the sets $\{x\in\Omega: d_A(x)=d_\Gamma(x)<r_j\}$ in the degeneracy condition.
This set  identifies the part of  the $r$-neighbourhood of $\Gamma$ which is closest to the set $A\subset F_j$.
Therefore $\delta_j$ is a bound on the degeneracy of the coefficients on the $j$-th face of the boundary.

The proof of Statement~\ref{tuni4.1-1} is by a slight elaboration of the proof of  Proposition~\ref{puni3.1}.

If $A$ is a bounded non-empty subset of $\Gamma$ and $A_j=A\cap \overline F_j$ then $A=\bigcup^m_{j=1}A_j$.
Now define $\eta_{r,n}$ as in the proof of Proposition~\ref{puni3.1} with $0<r\leq \inf_j r_j$.
Then
\[
\|\eta_{r,n}\|_2^2\leq |A_r|\leq \sum^m_{j=1}|A_{j,r}|\leq \sum^m_{j=1} c_j \,r^{d-s_j}
\]
for all small $r>0$ uniformly in $n$ where $A_{j,r}=\{x\in\Omega: d_{A_j}(x)<r\}$.
This follows by repetition of the reasoning in the proof of Proposition~\ref{puni3.1}.

Next set $\Lambda_{j,r}=\{x\in\Omega: d_{A_j}(x)=d_\Gamma(x)<r\}$.
Since $\supp\eta_{r,n}\subset\overline{A_r}$, one has 
\begin{eqnarray*}
h_N(\eta_{r,n})=\sum^m_{j=1} \int_{\Lambda_{j,r}}dx\,\Gamma(\eta_{r,n})(x)
&\leq& \sum^m_{j=1}b_j\int_{\Lambda_{j,r}\cap A_r}dx\, d_{F_j}(x)^{\delta_j}|(\nabla\eta_{r,n})(x)|^2\\
&\leq& \sum^m_{j=1}b_j\int_{A_{j,r}}dx\, d_{A_j}(x)^{\delta_j}|(\nabla\eta_{r,n})(x)|^2\\
&\leq&\,(\log n)^{-2}\sum^m_{j=1}b_j\int_\Omega dx\,\one_{\{x:r n^{-1}\leq d_{A_j}(x)\leq r\}}\,d_{A_j}(x)^{-(2-\delta_j)}
\;.
\end{eqnarray*}
Thus if $\delta_j\geq2$  for all $j$  one immediately has bounds 
\[
h_N(\eta_{r,n})\leq \sum^m_{j=1} b_j\,|A_{j,1}|\,(\log n)^{-2}
\]
uniformly for $r\leq1$. 
Then $\inf_{n\geq1}\|\eta_{r,n}\|_{D(h_N)}^2\leq \sum^m_{j=1} c_j \,r^{d-s_j}$ for all $r\leq 1$  and consequently $\capp_{h_0}(A)=0$.
Alternatively, if $\delta_j<2$ for one or more $j$ but $\delta_j\geq \delta_{c,j}$  then one deduces as in  the proof of Proposition~\ref{puni3.1}  
that $h_N(\eta_{r,n})\leq b\,(\log n)^{-1}$ for all $n\geq 2$ uniformly for $r\leq1$.
Again it follows that $\capp_{h_0}(A)=0$.
Since this is valid for all bounded subsets $A$ of $\Gamma$ it follows by monotonicity  that $\capp_{h_0}(\Gamma)=0$.
 Hence $h_D=h_N$ by Proposition~\ref{puni2.1}.

Statement~\ref{tuni4.1-2} follows directly from Proposition~\ref{puni3.4}.
One now chooses $z\in F_j$ and replaces $s$ by $s_j$ and $\delta$ by $\delta_j$.
\hfill$\Box$

\bigskip

Theorem~\ref{tuni1.1}  and the foregoing extensions  do not apply 
directly to domains where the boundary separates $\Omega$ locally.
More precisely, we define $\Gamma$ as separating $\Omega$ locally 
at $z\in \Gamma$ if there exist an $r>0$ and open subsets $\Omega^{1}_{z,r}$,
$\Omega^{2}_{z,r}$ of $\Omega$ such that 
$\Gamma_{\!z,r}$ separates $\Omega^{1}_{z,r}$ and
$\Omega^{2}_{z,r}$ and
$\Omega_{z,r}=\Omega^{1}_{z,r}\cup \Omega^{2}_{z,r}\cup \Gamma_{\!z,r}$.
Note that $\Omega^{1}_{z,r}$ and $\Omega^{2}_{z,r}$ need not be connected since
$\Omega_{z,r}$ can even contain  infinitely many components. 
If $\Gamma_{\!z,r}$ separates $\Omega$ locally then the Hausdorff dimension of 
$\Gamma_{\!z,r}$ is greater or equal to $d-1$.
Therefore if $\Gamma$ is Ahlfors $s$-regular one must have $s\geq d-1$.
The local separation rules out the possibility that  $\Omega_{z,r}$ is $\Omega$-uniform.
Nevertheless, it is in fact enough for Theorem~\ref{tuni1.1} that
there exists at least one $\Omega$-uniform component of $\Omega_{z,r}$,
say $\Omega^{1}_{z,r}$, but we omit further details.

\section{Illustrations and examples}\label{S5}

There are two constraints placed on the domain $\Omega$ in the foregoing discussion, Ahlfors regularity of the boundary and  local uniformity  near the boundary.
To conclude we describe general situations for which these properties are valid and  and specific examples for which they can fail.
We begin by examining the regularity property.

First assume  $\Omega$ is a Lipschitz domain.
Specifically assume that  each boundary section $A=\Gamma\cap B(z\,;R)$ with $z\in\Gamma$ and $R>0$
has a finite cover by balls $B_k=B(x_k\,;r_k)$, $k=1,\ldots, N$ such that each
subsection $A_k= A\cap B_k$ is, after a suitable rotation and translation of coordinates,  the graph of a Lipschitz continuous function $\varphi_k$.
Then $\Gamma$ has Hausdorff dimension $s= d-1$ and $\ch^s$  measures  the surface area of the boundary (see, for example, \cite{EvG} Chapters~2 and 3).
Moreover,  it follows by a standard calculation (see \cite{EvG}, 3.3.4B) that the Ahlfors regularity property (\ref{euni1.3}) is valid with $\mu=\ch^s$.
Since the boundary is locally `flat'  the $\Omega$-uniformity property is automatically satisfied for all $z\in\Gamma$.
Therefore  both the geometric assumptions on $\Omega$ in Theorem~\ref{tuni1.1} are satisfied with $s=d-1$.

Secondly,  consider a domain $\Omega$   whose boundary $\Gamma$ is 
 (locally) a  self-similar fractal (see \cite{Fal1} \cite{Fal}).
Specifically assume that  $F_1,\ldots, F_m$   are similarity transformations of a closed subset $D$ of $\Ri^d$, 
with similarity ratios $r_1,\ldots, r_m$, and that $\Gamma$ is the 
unique compact subset of $\Ri^d$ satisfying the self-similarity condition $\Gamma=\bigcup^m_{k=1}F_k(\Gamma)$.
If the $F_k$ satisfy the  open set condition introduced by Hutchinson \cite{Hut} then $\Gamma$ is Ahlfors $s$-regular with $s$  the unique solution of the
equation $r_1^s+\ldots+r_m^s=1$.
(For details see  \cite{Fal1}, Chapter~2, and especially Exercise~2.11, or \cite{Fal}, Chapter~9.)
These observations allow one  to construct a multitude of domains with Ahlfors regular boundaries but the local uniformity condition is not necessarily satisfied.
We will illustrate this with various standard examples.
The simplest and best known is perhaps the (modified) von Koch snowflake domain in two-dimensions.

\begin{figure}[!ht]
\begin{center}
\includegraphics[width=13cm]{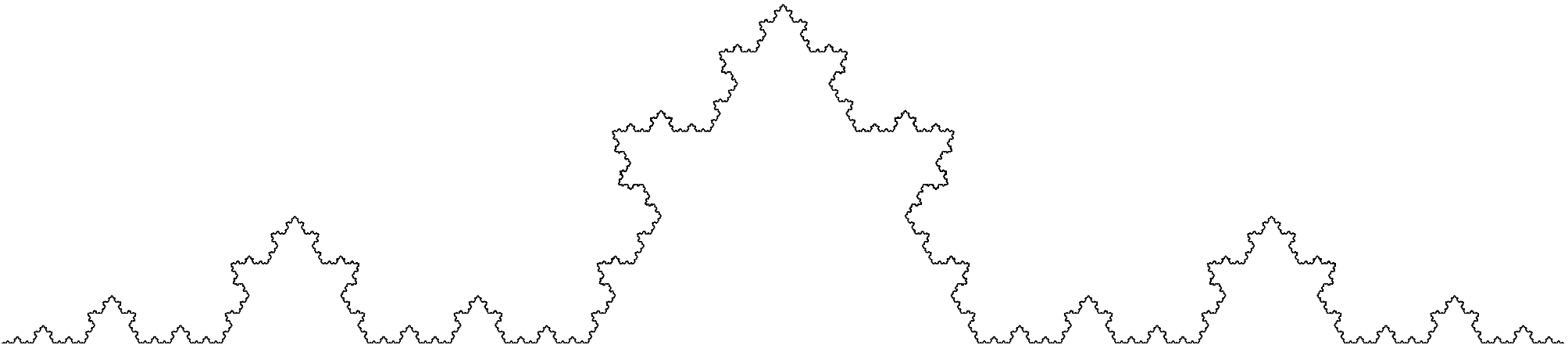}
\end{center}
\caption{Modified von Koch snowflake curve with $\lambda=1/4$.}
\end{figure}

\begin{exam}\label{exuni5.3}(von Koch snowflakes)
Let $\lambda\in\langle0,1/3]$ and let $E_0\subset \Ri^2$ be a  line segment of unit length.
Define $E_1$  by replacing the middle 
segment of $E_0$ with length $\lambda$ by the two sides of the equilateral triangle based
on the middle segment.
Thus $E_1$ is  the union of  four line segments  two of length $\lambda$ and two of length $(1-\lambda)/2$.
It is the image of $E_0$ under four similarity transformations with similarity ratios $(1-\lambda)/2,\,\lambda, \,\lambda, (1-\lambda)/2$.
Next define  $E_2$ by repeating this procedure on each of the four segments of $E_1$.
The (modified) von Koch `curve' corresponding to the choice of $\lambda$ is the self-similar set obtained by iteration of this 
procedure. (See \cite{Fal} Example~9.5 and Figure~1.)
Its Hausdorff dimension $s$ is the unique solution 
of $2\,\lambda^s+2((1-\lambda)/2)^s=1$.
Thus $s$ is a monotonically increasing function of $\lambda$ with values between $1$ and $\log4/\log3$.
In particular $s\to1$ as $\lambda\to 0$ and 
 $s=\log4/\log3$ if $\lambda=1/3$.
The  von Koch snowflake domain $\Omega$ is  the bounded,  simply connected, interior of the curve $\Gamma$
obtained by applying this construction to the  sides of a unit equilateral triangle with 
each von Koch curve pointing outwards.
(The standard construction is with $\lambda=1/3$.)
The boundary $\Gamma$  is Ahlfors $s$-regular
by the discussion preceding the example.
All  the von Koch domains and the corresponding unbounded exterior domains satisfy the  Martio-Sarvas uniformity
property.
An explicit proof of uniformity (for $\lambda=1/3$) is given in \cite{GyS}, Proposition~6.30.
Therefore Theorem~\ref{tuni1.1} applies directly to the quadratic form $h_0$ defined on the von Koch domains and since $d=2$
one concludes that $h_0$ is Markov unique if and only if the degeneracy parameter $\delta\geq s$.
\end{exam}

 Finally we consider two families of $d$-dimensional self-similar fractals which have very similar definitions but are of a quite different nature.
 Together these examples illustrate various possible features including the failure of the local uniformity property near the boundary.

\begin{exam}\label{exuni5.4}(Vicsek snowflakes)
Let $\lambda\in\langle0,1/2\rangle$ and let $E_0$ be the unit cube  centred at the origin of $\Ri^d$, i.e.\ 
$E_0=\{x\in\Ri^d:|x_j|\leq1/2, \;j=1,\ldots,d\}$. 
Define $E_1$ as the union of $2^d+1$ cubes consisting of the central cube of $E_0$ with edge length $1-2\lambda$  and the $2^d$  corner cubes 
of   edge length~$\lambda$.
Then $\Gamma\subset E_0$ is defined as the self-similar set obtained by iterating this procedure.
(The case $d=2$ and $\lambda=1/4$ is illustrated by Figure~0.5 in \cite{Fal}).
Thus $\Gamma$ is invariant under  $2^d+1$ similarity transformations, one with  similarity ratio $(1-2\lambda)$ and the remaining
$2^d$  with similarity ratio $\lambda$.
Consequently, the  Hausdorff dimension $s$ of $\Gamma$ is  determined by the equation $2^d\lambda^s+(1-2\lambda)^s=1$.
Therefore  $s\in\langle1,d\rangle$ and  the value of $s$ increases monotonically with $\lambda$.
Moreover,  $s\to1$ as $\lambda\to0$ and $s\to d$ as $s\to 1/2$.
Hence for each $d$ there is a critical value $\lambda_d$ for which $s<d-1$ if $\lambda<\lambda_d$ and $s\geq d-1$ if $\lambda\geq \lambda_d$.
The value of $\lambda_d$ increases monotonically with $d$ and $\lambda_d\to 1/2$ as $d\to\infty$, e.g. 
$\lambda_2=0$, $\lambda_3=1/3$, $\lambda_4=(\sqrt{21}-3)/4$.
Now let $\Omega$ be the complement of $\Gamma$ in $\Ri^d$.
The uniformity properties of $\Omega$ are quite different in the two cases $d=2$ and $d\geq 3$.

First if $d=2$ then $E_0$ is a unit square and $\Gamma$ consists of the two diagonals of the square `decorated' with diagonal `antennae' 
in a  pattern  repeated at smaller and smaller scales by the self-similar construction.
The antennae invalidate the uniformity property for $\Omega$ and the local uniformity property of Section~\ref{S1} fails at all points of $\Gamma$.
For example the sets $\Omega_{z,R}$ with $z\in\Gamma$ are either separated or split by the antennae.
Since $d=2$ and $s>1$ it follows that  Theorem~\ref{tuni1.1}  is not applicable in this case.

Secondly, if $d\geq 3$ then $\Gamma$ consists of the diagonals of the $d$-dimensional unit cube again decorated with antennae parallel to the diagonal directions
with smaller scale antennae.
It then follows that $\Omega$ is a uniform domain
because one can now choose paths between pairs of points in $\Omega$ which circumvent the antennae.
The uniformity can be verified by a  variation of the argument given in \cite{GyS}, Proposition~6.30, for the von Koch snowflake.
Therefore in this case Theorem~\ref{tuni1.1} is  applicable for all choices of $\lambda$.
\end{exam}

\begin{figure}[!ht]
\begin{center}
\includegraphics[width=10cm]{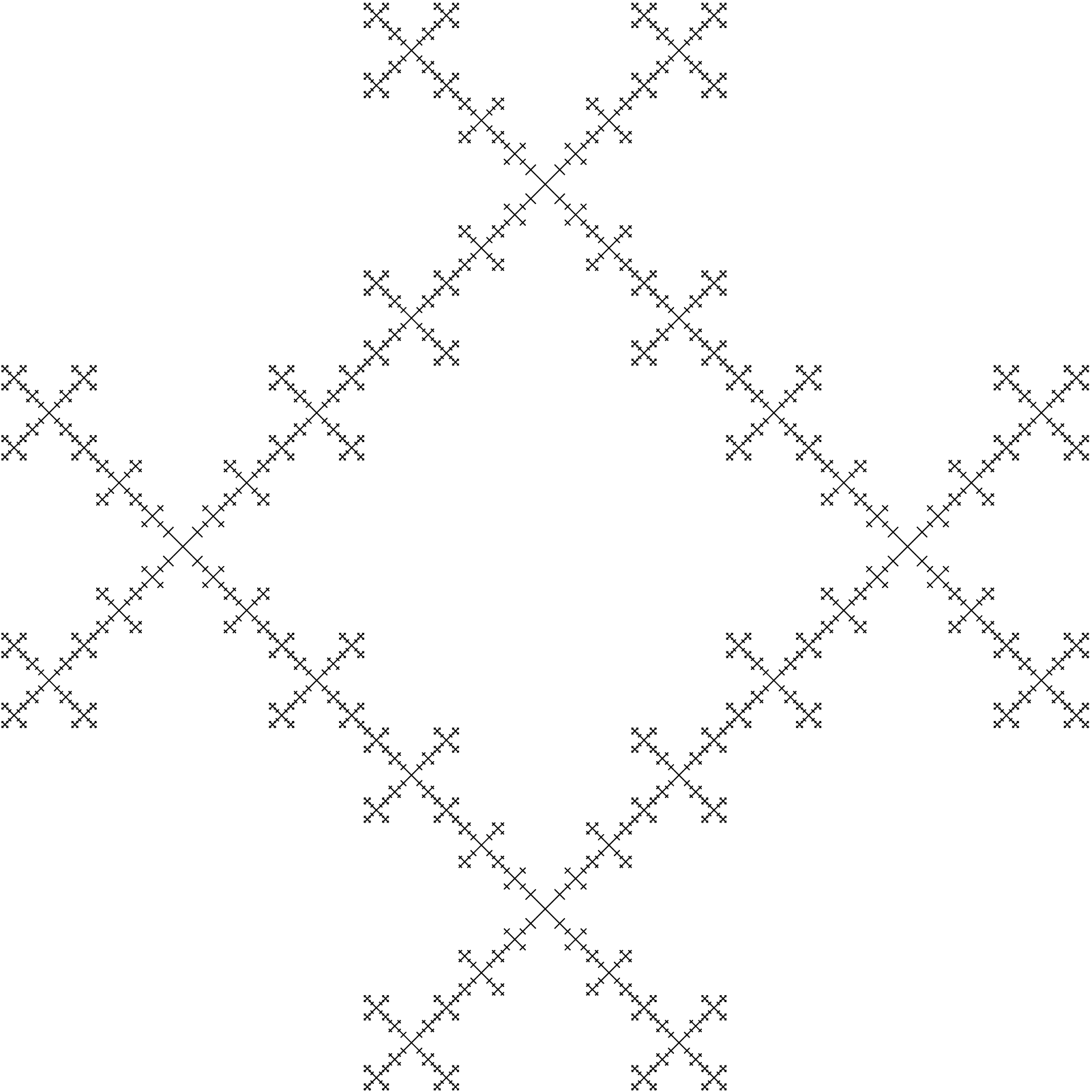}
\end{center}
\caption{Four copies of the Vicsek snowflake with $\lambda= 1/4$
bounding interior and exterior domains.}
\end{figure}

Note that if $d=2$ and one applies the Vicsek snowflake construction to four copies of the unit square $E_0$ centred at $(\pm1,0)$, $(0,\pm1)$, respectively, then the union of the  four copies of the self-similar set $\Gamma$  is a locally self-similar set
$\widetilde\Gamma$ which separates $\Ri^2$ into a bounded interior domain $\Omega_{\rm int}$ and an unbounded exterior domain $\Omega_{\rm ext}$ (see Figure~2). 
The common boundary $\widetilde\Gamma$  is Ahlfors $s$-regular but the uniformity properties fail for both domains.
If, however, $d\geq 3$ and one takes a $\Zi^d$-periodic partition of $\Ri^d$ into unit cubes   
and  then applies the  foregoing  construction to each cube the union of the resulting self-similar sets  connect at the corners of the cubes and form a periodic web $\widetilde\Gamma$.
This web  is a closed connected Ahlfors $s$-regular  set with $s\in\langle1,d\rangle$,  its complement  $\Omega=\Ri^d\backslash\widetilde \Gamma$ is a uniform  domain and 
 the conclusions of Theorem~\ref{tuni1.1} are valid for all possible values of $s$.

 Finally we consider the example of Cantor dust.
 The Cantor construction is similar to the Vicsek construction  but with quite different geometric properties.

 \begin{exam}\label{exuni5.6}(Cantor dust in $\Ri^d$)
Let $\lambda\in\langle0,1/2\rangle$ and  let $E_0$ be  a unit cube  in  $\Ri^d$.
Now define $E_1$ as the union of the $2^d$  corner cubes of $E_0$
with  edge length $\lambda$.
The self-similar set $\Gamma$  obtained by iterating this procedure
 is invariant under  $2^d$ similarity transformations with  similarity ratio $\lambda$.
 The set $\Gamma$ is a closed, completely disconnected, uncountable,  Ahlfors $s$-regular set with $s$ 
 determined by $2^d\lambda^s=1$.
 Thus the  Hausdorff dimension  is given by
 $d_H(\Gamma)=s=d\log2/\log(1/\lambda)$.
 Its value increases as $\lambda$ increases and takes all values in   $\langle0,d\rangle$.
Corollary~\ref{cuni1.1} now applies.
Thus $h_0$ is Markov unique if and only if $\delta\geq 1+(s-(d-1))$.
 
 One can also take a periodic partition of $\Ri^d$ into disjoint parallel cubes of side-length $2\lambda$ spaced at distance $1-2\lambda$
 and then apply the iteration procedure to each cube.
 The resulting set is a regularly spaced cloud of Cantor dust $\widetilde\Gamma$ with Hausdorff dimension $d\log2/\log(1/\lambda)$.
 Corollary~\ref{cuni1.1}  is again applicable.
 Therefore $h_0$ is Markov unique if and only if $\delta\geq 2+s-d$.
 \end{exam}

The last two examples are  interesting even for non-degenerate forms $h_0$, i.e.\ for the case $\delta=0$.
For example if $h_0(\varphi)=\sum^d_{j=1}\|\partial_j\varphi\|_2^2$ is the form of the Laplacian $\Delta$ and $\widetilde\Gamma$ is the Cantor dust cloud then
 $\Delta$ is Markov unique, or equivalently $L_1$-unique, on $\Ri^d\backslash\widetilde\Gamma$ if and only if 
$d\geq 2+s$.
Thus the presence of the cloud  influences the diffusion if and only if the roughness parameter $\lambda$ is sufficiently
large that $s>d-2$.
A similar conclusion is valid for the Vicsek web if $d\geq 3$.

\bigskip

\subsection*{Acknowledgement}
The authors are indebted to Adam Sikora who suggested 
that the Hardy--Rellich inequality could be used to identify the critical degeneracy for  Markov uniqueness.
The second named author is grateful to Michael and Louisa Barnsley for numerous discussions on fractal geometry.
%During the preparation of this work, the first named author was 
%partially supported by the Academy of Finland (grant no.\ 252108).
The first named author was 
partially supported by the Academy of Finland (grant no.\ 252108) during the preparation of this work.

%\bibliography{refbib}

\end{document}